\documentclass{article}
\usepackage{latexsym}
\usepackage[utf8]{inputenc}
\usepackage{amsmath,amssymb,amsthm,enumerate,bm,color}

\newtheorem{thm}{Theorem}[section]

\newtheorem{cor}[thm]{Corollary}
\newtheorem{prop}[thm]{Proposition}
\newtheorem{conj}[thm]{Conjecture}

\newtheorem{defn}[thm]{\bf Definition}

\theoremstyle{remark}
\newtheorem{rem}[thm]{\bf Remark}

\numberwithin{equation}{section}

\newcommand{\beq}{\begin{equation}}
\newcommand{\eeq}{\qed \end{equation}}

\def\1{\mathbf{1}}
\def\2{\mathbf{2}}

\begin{document}

\begin{Large}
\centerline{\bf The Zero Set of an Electric Field}
\centerline {\bf from a Finite Number of Point Charges}
\medskip

\date{August 12, 2026}

\centerline{\bf Tam\'as Erd\'elyi, Joseph Rosenblatt, Rebecca Rosenblatt}
\end{Large}
\medskip

\begin{abstract}  We consider the structure of the zero set in ${\mathbb R}^3$ of the electric vector field ${\mathbf F}=(X,Y,Z)$ from a finite set of point charges.  We are most interested in the case where the point charges all lie in a plane, and we consider just the zero set in ${\mathbb R}^2$ of the electric vector field ${\mathbf F}=(X,Y)$ from the finite set of point charges.  The conjecture is that the zero set of ${\mathbf F}=(X,Y)$ in ${\mathbb R}^2$ is finite.   We show fairly easily that this conjecture is true in a Special Case: when the point charges for ${\mathbf F}=(X,Y)$ lie on a line, and we consider the possible zeros throughout ${\mathbb R}^2$.  However, even in this Special Case, it is hard to get complete structural information about the zero sets of $X$ and $Y$ separately.  We describe structural information about the asymptotic directions at infinity of these two zero sets, and relate this to the interlacing of the zero sets of sequences of polynomials.  Then we consider the General Case where the point charges can be anywhere in the plane.  As in the Special Case, we construct sequences of polynomials whose zero sets include the asymptotic directions at infinity of the zero sets of $X$ and $Y$ separately.  But now the asymptotic directions are not necessarily interlacing, and the structure of the zero sets of these polynomials is less evident.  Nonetheless, using these polynomials it might be possible to show that the zeros of ${\mathbf F}=(X,Y)$ are bounded, and then perhaps also, as a result, confirm the conjecture that the zero set of ${\mathbf F}=(X,Y)$ is finite. 
\end{abstract}


\section{Introduction}\label{Intro} 

This article is about the zero sets of electric vector fields ${\mathbf F}=(X,Y,Z)$ in ${\mathbb R}^3$ from a finite set of point charges.  
While some general facts are clear in this setting, in order to get detailed proofs of various facts, we focus on the case where the charges are 
in a plane.  That is, ${\mathbf F}=(X,Y)$ is a finite point charge vector field in ${\mathbb R}^2$.  The well-known, still open conjecture is that in 
this case the zero set of ${\mathbf F}=(X,Y)$ is finite.  

Using the results in Yanushauskaus~\cite{Jan}, it is not difficult to prove this conjecture in the Special Case where the charges themselves lie on a 
line in ${\mathbb R}^2$, but the zeros are allowed to be anywhere in a bounded disk of the plane ${\mathbb R}^2$.  Thus, the Special Case may seem to be an ``easy" case to 
analyze completely.  But it turns out that this is not true if one wants to know more than just that the zero set of ${\mathbf F}=(X,Y)$ is finite in any bounded disk 
of ${\mathbb R}^2$. In particular, it is hard to give the detailed structure of the zero sets $X=0$ and $Y=0$ separately.  So besides showing that the zero set is finite in 
the Special Case, here we obtain structural information about the zero sets of $X$ and $Y$ at infinity by analyzing power series expansions, using binomial series, 
geometric series, and also critical moments that are given by the coefficients of the vector field and the positions of the charges.  Then we consider the General Case 
=where the point charges can be anywhere in the plane.  It is even more difficult to extend such detailed structural analysis of the zero sets $X=0$ and $Y=0$ separately 
to the General Case where the point charges are in general position in ${\mathbb R}^2$.  We present work toward showing that the zeros of ${\mathbf F}(X,Y)$ are bounded. 
If this boundedness is correct, we discuss what else would be needed to confirm the conjecture that the zero set of ${\mathbf F}=(X,Y)$ is finite.   

The outline of this article is as follows.  In this section we set up the notation and the problems we study in this paper. In Section~\ref{NewResults} 
we list most of our new results. In Section~\ref{Curves} we show that the 
zero set of ${\mathbf F}=(X,Y)$ is finite when the point charges lie on a line, that is in the Special Case.  In Section~\ref{Asymptotes}, we derive facts about the asymptotic directions at infinity of the zero sets of $X$ and $Y$ for the Special Case, and explain why most of these directions do not overlap.  Following this, in Section~\ref{GENERAL} we consider the case where the point charges can lie anywhere in the plane and develop formulas for sequences of polynomials whose zero sets contain the asymptotic directions at infinity of the zero sets of $X$ and $Y$.  Perhaps these can be used to show that the zero set of ${\mathbf F}=(X,Y)$ is bounded.  After that is known, a local argument is what would be needed to show that indeed the zero set of ${\mathbf F}=(X,Y)$ is finite.

\subsection{Finite Electric Fields in Three Dimensions}\label{GenThreeSpace}

\begin{defn}\label{FEF3}
Let $(x_j,y_j,z_j) \in {\mathbb R}^3,\,j=1,2,\ldots,M$, be distinct points with electric charges $a_j \in {\mathbb R}$ at $(x_j,y_j,z_j)$. The vector field ${\mathbf F}=(X,Y,Z)$ 
is called a finite point charge electric vector field defined on ${\mathbb R}^3$ by the component functions 
$$X(x,y,z) := \sum_{j=1}^M{\frac{a_j(x-x_j)}{r_j^3}}\,, \quad  Y(x,y,z) := \sum_{j=1}^M{\frac{a_j(y-y_j)}{r_j^3}}\,, \quad Z(x,y,z) := \sum_{j=1}^M\frac{a_j(z-z_j)}{r_j^3}\,,$$
with real numbers $a_j$, $j=1,2,\ldots,M$, which are not all $0$, and 
$$r_j := ((x-x_j)^2+ (y-y_j)^2 + (z-z_j)^2)^{1/2}\,, \qquad j=1,2, \dots,M\,.$$
\end{defn}

We know from Gibson~\cite{Gibson} that in ${\mathbb R}^3$, $X$, $Y$, and $Z$ can have zero sets consisting of points, curves, and even two-dimensional surfaces. 
However, how these intersect, which is where the full vector field is $(0,0,0)$, is not clear.  Are the sets of intersections also points, curves, and surfaces too, or just 
points and curves? Are the number of these intersections finite or infinite?  Also, how many points, curves, and surfaces will there be at most, say given only $M$?  

The basic structure is limited by Yanushauskas~\cite{Jan}. Yanushauskas results can be used because the electric field in the three-dimensional case is the gradient of a 
harmonic function. This shows that the common zero set consists of a locally finite set of points and analytic curves.  That is, even though a component can be zero on a 
two-dimensional surface, the full vector field cannot be $(0,0,0)$ on a two-dimensional surface.  

\begin{defn}  We use the terms curves and surfaces to mean at least ones that can be continuously parameterized. But actually, the curves and surfaces that occur in this 
article are even more special because not only are they analytic, but also they are branches in algebraic varieties.
\end{defn}

From Yanushauskas~\cite{Jan}, we have the following.

\begin{prop}\label{Yanush} A finite electric vector field ${\mathbf F}=(X,Y,Z)$ on ${\mathbb R}^3$ does not vanishes identically on a set containing a 
two-dimensional surface.
\end{prop}

But we still have structural issues to consider in order to understand the geometry of such a zero set.
\medskip

\noindent {\bf Structure Conjecture}: The zero set of $X$, $Y$, or $Z$ is asymptotically planar or linear at infinity, with only a finite number of planes and lines 
being possible.  Moreover, these planes or lines are distinct among the three components and the overlap of these three zero sets can only consist of a 
finite set of points and asymptotically linear curves.
\medskip

Of course, we could have a zero set containing asymptotically a surface or curve in one direction, and another surface or curve in another direction.  
But if the Structure Conjecture is true then the common zero set would either be bounded or consist of a finite number of curves each of which is asymptotically 
linear at infinity.  

\begin{rem} Examples show that the zero set in the three-dimensional case can contain curves (e.g. circles and straight lines).  Counting these and understanding 
what types of curves one can get may be a future project.  The main question otherwise is whether or not the unbounded portion is asymptotically linear.  
\end{rem}

\begin{rem} With some additional work we can somewhat address the Structure Conjecture by using the basic theory of varieties.  
We could use Whitney~\cite{WI}, Whitney~\cite{WII}, and Whitney and Bruhat~\cite{WB}.  But this will only give us information about the zero sets $X=0$, $Y=0$, and $Z=0$ 
being contained locally in nonidentically zero real analytic varieties.
Instead, we will use a product argument, and the structure of real algebraic varieties, to get more information. See Corollary~\ref{inside} and Whitney~\cite{WI} and Gibson~\cite{Gibson} 
for background. 
\end{rem}

\begin{rem}  It is an interesting issue that deriving the polynomial $P$ in Proposition~\ref{squaring} is not so straightforward if one wants to use a different method than the one used
in the prof of Proposition~\ref{squaring}. Most elementary cases of equations with square roots are handled just by isolating the terms on one side of the equality, squaring, rearranging 
the terms, and continuing this process.  But with more than three terms this does not proceed so easily.  See Killian~\cite{K} for an explanation of how to make the method of multiple 
squaring work to show that the zero set is a subset of a nonidentically zero real algebraic variety.
\end{rem}

\medskip

\begin{rem} It is possible that in fact in Proposition~\ref{inside} not only is the zero set contained in the associated structure, but it is itself in fact of the same structure.  
But it is possible that we might have the zero sets, or the component zero sets, containing an infinite, discrete set of points which are all contained in some algebraic curve. 
Computations and the results that we have suggest that this is never the case, but we do not yet have a rigorous argument that shows this.
\end{rem}

\subsection{Finite Electric Fields, General Case in the Plane}\label{GenPlane}

\begin{defn}\label{FEF2}
Let $(x_j,y_j) \in {\mathbb R}^2, j=1,2,\ldots,M$, be distinct points with charges $a_j \in {\mathbb R}$ at $(x_j,y_j)$. The vector field ${\mathbf F}=(X,Y)$
is called a finite point charge electric vector field defined on ${\mathbb R}^2$ by the component functions
$$X(x,y) := \sum_{j=1}^M{\frac{a_j(x-x_j)}{r_j^3}}\,, \qquad  Y(x,y) := \sum_{j=1}^M{\frac{a_j(y-y_j)}{r_j^3}}\,,$$
with real numbers $a_j, j=1,2,\ldots,M$, which are not all $0$, and
$$r_j := ((x-x_j)^2+ (y-y_j)^2)^{1/2}\,, \qquad j=1,2, \dots,M\,.$$
\end{defn}

Examples suggest that now the zero sets of $X$ and $Y$ consists of points and curves, both bounded and unbounded. However, the common zero set is more limited.
\medskip

\noindent {\bf Conjecture}\label{MAINCONJ}: The zero set $X=Y=0$ of a two-dimensional finite electric vector field ${\mathbf F}=(X,Y)$ consists of a finite set of points.
\medskip

\noindent Of course, we would want to know also what the the answer to the corresponding Maxwell's question  would be: how many points can be in the zero set 
$X=Y=0$ of an electric vector field ${\mathbf F}(X,Y)=(X,Y)$, given only $M$?  
See Gabrielov, Novikov, and Shapiro~\cite{GNS}.  See also Proposition 4.1 in \cite{AHKT} and Killian~\cite{K}.

Despite a number of attempts to find simple method to resolve {\bf Conjecture}~\ref{MAINCONJ} and to estimate the number of points, the conjecture has not been proven.  

One can consider two separate issues with describing the zero set of a finite electric vector field ${\mathbf F}=(X,Y)$.  One issue is how to show that the zero set of 
${\mathbf F}=(X,Y)$ must be in a large disk. This is not the case for the zero sets $X=0$ and $Y=0$ separately; they typically do have portions arbitrarily far from the origin. 

The second issue is how to show that the zero set of a finite electric vector field ${\mathbf F}=(X,Y)$ is finite. The second issue is not completely resolved 
at this time unless we had a result like Proposition 4.1 in \cite{AHKT}.

Actually, what we seek to prove, with some technical work that is much more difficult to complete than would be desirable, 
is that at infinity the zero sets $X=0$ and  $Y=0$ consist of a finite number of curves that are to some degree at least 
asymptotically linear. In the process, or as a result, we would want to find equations that give the slopes of these asymptotic 
lines and show that they are distinct from one another when we consider $X$ in contrast to when we consider $Y$. This is how we could 
show that the common zeros of $X$ and $Y$ must lie in a large disk.
  
The obvious approach to these questions would seem to be just a smart use of implicit differentiation.  For example, where $X=0$, if there is a solution $y = y(x)$ 
or $x = x(y)$, we should be able to compute $dy/dx$ or $dx/dy$ and work with this.  In particular, if we believe that a given curve $y = y(x)$ in the zero set $X=0$ 
must be asymptotically linear at infinity, we would just need to show that $y$ is asymptotically $\alpha x$, for some $\alpha$.  At this time, we cannot see how to carry 
this out in general.

So instead, what we do is consider the implicit equation $X=0$ and expand it as an infinite series.  This series has to be adapted to the region in question; 
one series does not seem to suffice for the whole plane.  Then we consider these series asymptotically and get equations that determine the asymptotic directions of the 
zero set $X=0$ at large distances from the origin.  These equations only have a finite number of solutions.  Then we carry out the same thing for the zero set $Y=0$.  
We then may be able to show generally that the asymptotic directions are different for the two zero sets.  This approach gives us at least part of the results that we want.   
In our study of finite electric vector fields on the plane in Section~\ref{GENERAL} the following definition of the number $L$ is crucial.

\begin{defn}\label{LGen} Suppose $a_j, j=1,2,\ldots,M$, are real numbers which are not all zero. Let $0 \geq L$ be the largest integer for which
$$\sum_{j=1}^M{a_jx_j^ly_j^i} = 0, \qquad i+l < L\,.$$
\end{defn}

\noindent Such an integer $0 \leq L \leq 2M-2$ exists, otherwise
\begin{equation} \notag \sum_{j=1}^M{a_jr_j^{2k}} = \sum_{j=1}^M{a_j((x-x_j)^2 + (y-y_j)^2)^k} = \sum_{j=1}^M{a_j(x^2 - 2x_jx + x_j^2 + y^2 - 2y_jy + y_j^2)^k} = 0
\end{equation}
for all $k=0,1,\ldots,M-1$ and $(x,y) \in {\mathbb R}^2$, and by choosing a point $(x,y)$ such that the values $r_j$, $j=1,2,\ldots,M$, are distinct,
the non-vanishing property of the Vandermonde determinants would imply that $a_1=a_2=\cdots=a_M=0$.

\subsection{A Special Case in the Plane}\label{Special}

Let us clarify the language that will follow below.  The {\em Special Case} will mean that we take points in ${\mathbb R}^2$ that are all on a line and consider the zeros of ${\mathbf F}=(X,Y)$.  

\begin{defn}\label{GC}
The {\em General Case} will mean that we take any distinct points $(x_j,y_j) \in {\mathbb R}^2$ and consider the zeros of ${\mathbf F}=(X,Y)$. 
\end{defn}

By rotating and translating, or by an appropriate choice of the coordinate system, the Special Case becomes the following.  

\begin{defn}\label{SC}{Special Case} Let $(x_j,0), j = 1,2,\ldots,M$, be distinct points on the positive $x$-axis. Assume that $0 < x_1 < x_2 < \cdots < x_M$. 
Put electric charges $a_j \in {\mathbb R}$ at the points $(x_j,0)$. Then the vector field ${\mathbf F}(X,Y)$ on the plane ${\mathbb R}^2$ is given by 
$$X(x,y) := \sum_{j=1}^M{\frac{a_j(x-x_j)}{\left( (x-x_j)^2 + y^2 \right)^{3/2}}}, \qquad  Y(x,y) = \sum_{j=1}^M{\frac{a_jy}{\left( (x-x_j)^2 + y^2 \right)^{3/2}}}\,.$$  
\end{defn}

Where is ${\mathbf F}=(X,Y)=(0,0)$? 
Can we prove that $X$ and $Y$ only simultaneously at only a finite number of values $(x,y) \in {\mathbb R}^2$?

One aspect of this analysis that is worth noting is that certain {\em moments} are closely tied to the structure of the zero sets.

\begin{defn}\label{intromoments}
In the Special Case we call $\mu_k:=\sum\limits_{j=1}^M{a_jx_j^k}$ the $k$th {\em moment}, $k=0,1,\ldots\,.$. 
\end{defn}

If $k=0$, then we get the ground state moment $\mu_0 = \sum\limits_{j=1}^M{a_j}$. If $M$ is large,
then many of these moments could be zero without forcing all the charges $a_j$ to be zero.  Given an electric vector field in the Special Case, the critical $L$th moment 
is defined as follows. 

\begin{defn}\label{L} Given an electric vector field in the Special Case, there is a unique integer $0 \leq L \leq M-1$  such that 
$$\mu_k = \sum\limits_{j=1}^M{a_jx_j^k} = 0\,, \qquad k=0,1,\ldots,L-1\,, \qquad \text{and} \qquad \mu_L = \sum\limits_{j=1}^L{a_jx_j^L} \neq 0\,.$$ 
\end{defn}

\noindent The existence of such a unique number $0 \leq L \leq M-1$ follows from the non-vanishing property of the Vandermonde determinants. Indeed, if 
$\mu_k=0$ for all $j=0,1,\ldots,M-1$, then we are lead to an $M$ by $M$ system of linear equations 
$$\mu_k = \sum\limits_{j=1}^M{x_j^ka_j} = 0\,, \qquad k=0,1,\ldots,M-1\,,$$  
with an $M$ by $M$ matrix formed by the coefficients $x_j^k \, j=0,1,\ldots,M-1, \, k=0,1,\ldots,M-1$. As $0 < x_1 < x_2 < \cdots < x_M$, the 
determinant of this matrix is a non-vanishing Vandermonde determinant. This implies that it has only the trivial solution $a_1=a_2=\ldots=a_M=0$, 
which contradicts to our assumption that not all $a_j, \, j=1,2,\ldots,M$, are $0$.

As the number of point charges $M$ and the number $0 \leq L \leq M-1$ increase, the structure of the zero sets of $X$ and $Y$ become more complicated,
Without carrying out an asymptotic analysis using power series expansions,
we have not found a way to see that these numbers $L$ play an important role in the structure of the zero sets $X$ and $Y$.

Given an electric vector field in the Special Case, we use the choice of $L$ given by Definition~\ref{L} throughout Sections~\ref{SecYTypeI}, ~\ref{SecXTypeI}, ~\ref{SecYTypeII}, 
and ~\ref{SecXTypeII}. Typically, as $L$ gets larger, the complexity to detect possible asymptotic directions for the zero sets $X=0$ and $Y=0$ increases.

\section{New Results}\label{NewResults}

We start with the following simple observation.

\begin{prop} 
In the Special Case the function ${\mathbf F}=(X,0)$ vanishes only at finitely many $x \in {\mathbb R}$.
\end{prop}  

\begin{proof}
This is the subcase that $y = 0$ of the Special Case. This is actually the case of a finite electric vector field on a line.  But now $Y(x,y) = Y(x,0) = 0$ and 
$$X(x,y) = X(x,0) = \sum_{j=1}^M{\frac{a_j(x-x_j)}{\left( (x-x_j)^2 \right)^{3/2}}} = \sum_{j=1}^M{\frac{a_j(x-x_j)}{|x-x_j|} \frac{1}{(x-x_j)^2}}\,.$$
Let $I_0 = (-\infty,x_1)$, $I_M = (x_M,\infty)$, and $I_j = (x_j,x_{j+1}), \enskip j=1,2,\ldots,M-1$.
The function $h(x) = \displaystyle{\frac{x-x_j}{|x-x_j|}}$ is identically $-1$ or identically $1$ on any of these intervals $I_j$ for $j=0,1,\ldots,M$.
Hence, $X(x,y) = X(x,0)$ is a rational function on any of these $I_j$ and can only have a finite number of zeros on $I_j$. This means that $F(x,0)$
can be zero only finitely many times too.
\end{proof}

In Section~\ref{subsec3.1} we prove the following.

\begin{prop}\label{XXX} In the Special Case the zero set of ${\mathbf F}=(X,Y)$ in every bounded disk of ${\mathbb R}^2$ is finite.
\end{prop}

In Section~\ref{subsec3.2} we outline the proof of the following two results. 

\begin{prop}\label{FINITE} In the Special Case the zero set of ${\mathbf F}=(X,Y)$ in ${\mathbb R}^2$ is finite.
\end{prop}

\begin{prop}\label{UPPERBOUND} In the Special Case the zero set of ${\mathbf F}=(X,Y)$ contains at most $9M^24^M$ points.
\end{prop}

In Sections ~\ref{subsec3.3} and ~\ref{subsec3.4} we give a detailed of the Proposition~\ref{FINITE}, while in Section ~\ref{subsec3.5} we give a detailed proof 
of Proposition~\ref{UPPERBOUND}.

In Section~\ref{SecYTypeI} we prove 

\begin{prop}\label{YTypeI} Let $\delta \in (0,1/2)$ be fixed in the Special Case. If
$$Y(p_m,q_m) = 0, \quad \left| \frac{p_m}{q_m} \right| < 1-\delta, \quad q_m \neq 0, \quad |q_m| \rightarrow \infty, \quad \lim_{m \rightarrow \infty}{\frac{p_m}{q_m}} = \beta\,,$$
then
$$\frac{d^L}{d z^L} \left( \frac{1}{(1+z^2)^{3/2}} \right)\bigg |_\beta = 0\,,$$
where $0 \leq L \leq M-1$ is defined by Definition~\ref{L}.
\end{prop}

In Section~\ref{SecXTypeI} we prove

\begin{prop}\label{XTypeI} Let $\delta \in (0,1/2)$ be fixed in the Special Case If
$$X(p_m,q_m) = 0, \quad \left| \frac{p_m}{q_m} \right| < 1-\delta, \quad q_m \neq 0, \quad |q_m| \rightarrow \infty, \quad \lim_{m \rightarrow \infty}{\frac{p_m}{q_m}} = \beta\,,$$
then
$$\frac{d^L}{d z^L} \left( \frac{z}{(1+z^2)^{3/2}} \right) \bigg |_\beta= 0\,,$$
where $0 \leq L \leq M-1$ is defined by Definition~\ref{L}.
\end{prop}

In Section~\ref{SecYTypeII} we prove

\begin{prop}\label{YTypeII}  Let $\delta \in (0,1/4)$ be fixed in the Special Case. If
$$Y(p_m,q_m) = 0, \quad |p_m| > x_M, \quad q_m \neq 0, \quad \left| \frac{q_m}{p_m} \right| < 1-2\delta, \quad |p_m| \rightarrow \infty, \quad 
\lim_{m \rightarrow \infty}{\frac{q_m}{p_m}} = \alpha\,,$$
then
$$\frac{d^L}{dz^L} \left( \frac{z^{L+2}}{(1+z^2)^{3/2}} \right)\bigg|_\alpha = 0\,,$$
where $0 \leq L \leq M-1$ is defined by Definition~\ref{L}.
\end{prop}

In Section~\ref{SecXTypeII} we prove

\begin{prop}\label{XTypeII} Let $\delta \in (0,1/4)$ be fixed in the Special Case. If
$$X(p_m,q_m) = 0, \quad |p_m| > x_M, \quad q_m \neq 0, \quad \left| \frac{q_m}{p_m} \right| < 1-2\delta, \quad |p_m| \rightarrow \infty, \quad \lim_{m \rightarrow \infty}{\frac{q_m}{p_m}} = \alpha\,,$$
then 
$$\frac{d^L}{dz^L} \left( \frac{z^{L+1}}{(1+z^2)^{3/2}} \right)\bigg |_\alpha = 0\,,$$
where $L$ is defined by Definition~\ref{L}.
\end{prop}

We think that the results in Sections~\ref{lace} and ~\ref{invlace} are interesting for their own merit. 
Using these results, in Section~\ref{SpecialSum} we can easily prove 

\begin{prop}\label{DISTINCT} In the Special Case the set of possible asymptotic directions of the zero set $X=0$ and the set of possible asymptotic directions 
of the zero set $Y = 0$ are distinct, except possibly the directions given by the lines $y = \pm x$.
\end{prop}

Let ${\mathbb R}[x_1,x_2,\ldots,x_n]$ denote the ring of polynomials in the variables $x_j$ with real coefficients.
The next result, Proposition~\ref{squaring} can be proved using some basic Galois theory, but in subsec~/ref{squaring} we give a proof that requires elementary theoretical discussions.
We used this same idea in a somewhat different way in the proof of Proposition~\ref{FINITE}.

\begin{prop}\label{squaring} Let
$$R = \sum_{j=1}^M{\frac{A_j}{B_j^{1/2}}}\,,  \qquad A_j,B_j \in {\mathbb R}[x_1,x_2,\ldots,x_n]\,,$$
be nonidentically zero on ${mathbb R}^n$, where $B_j \ge 0$ on ${\mathbb R}^n$ for all $j=1,2,\dots,M$.  Then the zero set
$R=0$ is a subset of a nonidentically zero real algebraic variety in ${\mathbb R}^n$.
\end{prop}

Combining Corrolary~\ref{squaring} and Gibson~\cite{Gibson} we get the following. 

\begin{cor}\label{inside} The zero sets of the components $X$ and $Y$ of a finite point electric field in ${\mathbb R}^2$ are contained in the union of a finite 
set of points and a finite set of algebraic curves in ${\mathbb R}^2$. The zero sets of the components $X$, $Y$, and $Z$ of a finite point electric field in ${\mathbb R}^3$ 
are contained in the union of a finite set of points, a finite set of algebraic curves, and a finite set of algebraic surfaces in ${\mathbb R}^3$.
\end{cor}

In Section~\ref{GENERAL} we focus on a study the structure of the zero set of a nontrivial finite point charge electric vector field ${\mathbf F}=(X,Y)$ in
${\mathbb R}^2$ for which the $M$ point charges can be in a general position. This section can be followed independently of the rest of the paper. We establish 
equations satisfied by the possible directions for the zero sets of $X$ and $Y$ separately, and we show that there are only finitely many possible 
asymptotic directions for both of these zero sets, given again by the zero sets of sequences of certain polynomials.
One result that these computations give, which does not seem to be obvious without the computations, is that the component zero sets have a finite number of 
asymptotic directions with a bound on the number of directions that is $4M-1$. We do not list the new results achieved in Section~\ref{GENERAL} here.  

\section{Finiteness of the Zero Set of ${\mathbf F}=(X,Y)$ in the Special Case}\label{Curves}


\subsection{Proof of Proposition~\ref{XXX}\label{subsec3.1}}

\begin{proof}  Suppose this is not the case.  Then there would be distinct, common zeros $(u_m,v_m)$ converging to a common zero $(u_0,v_0)$. 
Consider our electric field extended to ${\mathbb R}^3$ by taking the position of the point charges to be $(x_j,0,0)$ instead of $(x_j,0)$. Then $(x_m,y_m,0)$
are distinct, common zeros of this extension converging to the common zero $(x_0,y_0,0)$ of this extension. Let $C_m$ be the circle
that is obtained by rotating $(x_m,y_m, 0)$  about the x-axis (the axis where the point charges are position).  By rotational
symmetry, then all $C_m, m \geq 1$, lie in the common zero set in ${\mathbb R}^3$.  This contradicts Theorem 2 in Yanushauskas~\cite{Jan} which states
the fact that this common zero set is locally a finite set of points and a finite number of pieces of analytic curves.
\end{proof}

\subsection{Outline of the Proof of Proposition~\ref{FINITE}}\label{subsec3.2}

For the sake of completeness we present one of our self-contained proofs showing that in the Special Case there cannot be any non-trivial curves
in the zero set ${\mathbf F}=(X,Y)=(0,0)$. First, in Section~\ref{2Curves}, we prove that this zero set is, in fact, countable. The advantage of this 
argument is that we do not even need to define what we mean by ``non-trivial curve". Then we show that the zero set ${\mathbf F}=(X,Y) = (0,0)$ is a subset 
of an algebraic variety not containing a non-trivial curve. Such an algebraic variety contains only finitely many points. See Gibson~\cite{Gibson}.
Then, using B\'ezout's Theorem, (see Gibson~\cite{Gibson}), we give an upper bound for the number of points in the zero set ${\mathbf F}=(X,Y)=(0,0)$.

\subsection{The Zero Set ${\mathbf F}=(X,Y)$ Is Countable}\label{2Curves}\label{subsec3.3}

With $0 < x_1 < x_2 < \cdots < x_M\,$, let
$$X_m(x,y) := \sum_{j=1}^M{\frac{a_j(x-x_j)}{((x-x_j)^2 + y^2)^{(2m+1)/2}}}\,, \qquad m=1,2, \ldots\,,$$
and
$$Y_m(x,y) := \sum_{j=1}^M{\frac{a_jy}{((x-x_j)^2 + y^2)^{(2m+1)/2}}}\,, \qquad m=1,2, \ldots\,.$$
We have
\begin{equation} \notag \begin{split}
\frac{\partial{X_m}}{\partial x} & = \, \sum_{j=1}^M{a_j \, \frac{((x-x_j)^2 + y^2)-(2m+1)(x-x_j)^2}{((x-x_j)^2 + y^2)^{(2m+3)/2}}}\,, \cr
\frac{\partial{X_m}}{\partial y} & = \, \sum_{j=1}^M{a_j \, \frac{-(2m+1)y(x-x_j)}{((x-x_j)^2 + y^2)^{(2m+3)/2}}}\,, \cr
\frac{\partial{Y_m}}{\partial x} & = \, \sum_{j=1}^M{a_j \, \frac{-(2m+1)y(x-x_j)}{((x-x_j)^2 + y^2)^{(2m+3)/2}}}\,, \cr
\frac{\partial{Y_m}}{\partial y} & = \, \sum_{j=1}^M{a_j \, \frac{((x-x_j)^2 + y^2)-(2m+1)y^2}{((x-x_j)^2 + y^2)^{(2m+3)/2}}}\,. \cr
\end{split} \end{equation}
Let
$$H := \{(x,y): \enskip X_1=Y_1=0, \enskip y \neq 0\}$$
and
$$H_m := \{(x,y): X_u = Y_u = 0, \enskip y \neq 0 \enskip u=1,2,\ldots,m, \enskip X_{m+1}^2 + Y_{m+1}^2 > 0\}\,,$$
$m=1,2,\ldots$. Observe that for $(x,y) \in H_m$ we have
$$\frac{\partial{X_m}}{\partial x} = -\frac{\partial{Y_m}}{\partial y} = (2m+1)yY_{m+1}$$
and
$$\frac{\partial{X_m}}{\partial y} = \frac{\partial{Y_m}}{\partial x} =-(2m+1)yX_{m+1}\,.$$
In this section we use only undergraduate real analysis such as Vandermonde determinants and the Implicit Function Theorem and
prove a structural property of the set $H$. In particular, we prove that $H$ is a countable set. The basic idea is under the assumptions
$$X_1(x_0,y_0) = Y_1(x_0,y_0)=0\,, \qquad \frac{\partial X_1}{\partial x}(x_0,y_0) \neq 0$$
the algebraic varieties
$$\{(x,y):X_1(x,y) = 0\} \qquad \text{and} \qquad \{(x,y):Y_1(x,y)=0\}\,,$$
which are graphs of functions in a neighborhood of $(x_0,y_0)$ by the Implicit Function Theorem, intersect
each other in an orthogonal fashion. However, proving our results in this section requires more technical
details.

\begin{prop}\label{union}
We have
$$H = \bigcup_{m=1}^{M-1}{H_m}$$
unless $a_1=a_2= \cdots =a_M = 0$.
\end{prop}

\begin{proof}
Suppose $(x_0,y_0) \notin H$. For the sake of brevity let
$$r_j := ((x_0-x_j)^2 + y_0^2)^{1/2}\,, \qquad j=1,2,\ldots,M\,.$$
We have
\begin{equation}\label{3.1}
\sum_{j=1}^M{\frac {a_j}{r_j} \frac 1{(r_j^2)^k}} = 0\,, \qquad k = 1,2,\dots,M\,,
\end{equation}
and
\begin{equation}\label{3.2}\sum_{j=1}^M{\frac {a_jx_j}{r_j} \frac 1{(r_j^2)^k}} = 0\,, \qquad k=1,2,\dots,M\,.
\end{equation}
Let $\{R_1 < R_2 < \ldots < R_m\} = \{r_1,r_2,\ldots,r_M\}$.
Since our point charges $(x_j,0)$ are on a line, for every $u \in \{1,2,\ldots,m\}$ there are at most $2$ values of
$j \in \{1,2,\ldots,M\}$ for which $r_j = R_u$. Hence by ~(\ref{3.1}) and ~(\ref{3.2}) we have
\begin{equation}\label{3.3}
\sum_{u=1}^m{\frac{b_u}{R_u} \frac 1{(R_u^2)^k}} = 0\,, \qquad k = 1,2,\dots,m\,,
\end{equation}
and
\begin{equation}\label{3.4}
\sum_{u=1}^m{\frac {c_u}{R_u} \frac 1{(R_u^2)^k}} = 0\,, \qquad k=1,2,\dots,m\,,
\end{equation}
where
\begin{equation}\label{3.5} b_u = a_j \qquad \text{and} \qquad c_u = x_ja_j
\end{equation}
if there is only one value $j \in \{1,2,\ldots,M\}$ such that $r_j = R_u$, and
\begin{equation}\label{3.6}
b_u = a_{j_1} + a_{j_2} \qquad \text{and} \qquad c_u = x_{j_1}a_{j_1} + x_{j_2}a_{j_2}
\end{equation}
if there are two distinct values $j_1, j_2 \in \{1,2,\ldots,M\}$ such that $r_{j_1} = r_{j_2} = R_u$.
Using ~(\ref{3.3})--~(\ref{3.6}) and the well known non-vanishing property of Vandermonde determinants we conclude that
$$b_u = c_u = 0\,, \qquad u=1,2,\ldots,m\,,$$
from which
$$a_j = 0, \qquad j=1,2,\ldots,M\,,$$
follows.
\end{proof}

\begin{prop}\label{countable}
Each point in $H_m$ is an isolated point in $H_m$ for every $m=1,2,\ldots,M-1$, unless $a_1=a_2= \cdots =a_M = 0$.
\end{prop}

\begin{proof}
Let $m$ be a fixed positive integer. Suppose $(x_0,y_0) \in H_m$, that is,
either
\begin{equation}\label{3.7} X_u(x_0,y_0) = Y_u(x_0,y_0) = 0\,, \enskip y_0 \neq 0\,, \enskip u=1,2,\ldots,m\,, \enskip X_{m+1}(x_0,y_0) \neq 0\,,
\end{equation}
or
\begin{equation}\label{3.8} X_u(x_0,y_0) = Y_u(x_0,y_0) = 0\,, \enskip y_0 \neq 0\,, \enskip u=1,2,\ldots,m\,, \enskip Y_{m+1}(x_0,y_0) \neq 0\,.
\end{equation}

If ~(\ref{3.7}) holds, then we have
\begin{equation}\label{3.9} \frac{\partial{X_m}}{\partial y}(x_0,y_0) = -(2m+1)y_0X_{m+1}(x_0,y_0) \neq 0
\end{equation}
and
\begin{equation}\label{3.10} \frac{\partial{Y_m}}{\partial x}(x_0,y_0) = -(2m+1)y_0X_{m+1}(x_0,y_0) \neq 0
\end{equation}  
Also,
\begin{equation}\label{3.11} \begin{split} \frac{\partial{Y_m}}{\partial y}(x_0,y_0) & = \, \frac{Y_m(x_0,y_0)}{y_0} - (2m+1)y_0Y_{m+1}(x_0,y_0) \cr
& = \, -(2m+1)y_0Y_{m+1}(x_0,y_0) \cr \end{split}
\end{equation}
and
\begin{equation}\label{3.12} \begin{split} \frac{\partial{X_m}}{\partial x}(x_0,y_0) & = \, (2m+1)y_0Y_{m+1}(x_0,y_0) - \frac{(2m)Y_m(x_0,y_0)}{y_0} \cr
& = \, (2m+1)y_0Y_{m+1}(x_0,y_0)\,. \end{split}
\end{equation}

By the Implicit Function Theorem and ~(\ref{3.9}), there are $\delta > 0$, $\varepsilon > 0$, and a function $y = f(x)$ differentiable on
$(x_0-\delta,x_0+\delta)$ such that the set
\begin{equation} \notag \begin{split} C_1 : & = \, \{(x,y): X_m(x,y) = 0, \enskip x \in (x_0-\delta,x_0+\delta), \enskip y \in (y_0-\varepsilon,y_0+\varepsilon)\} \cr
& = \, \{(x,f(x)): \enskip x \in (x_0-\delta,x_0+\delta)\}\cr \end{split}
\end{equation}
is the graph of a differentiable function $y = f(x)$ on $(x_0-\delta,x_0+\delta)$, and the equation for the tangent line to the curve $C_1$ at the point $(x_0,y_0)$ is
$$y-y_0 = \alpha(x-x_0)\,,$$
where
$$ \alpha = \left( \frac{\partial{X_m}}{\partial x} \Big / \frac{\partial{X_m}}{{\partial y}} \right)(x_0,y_0
= \frac{(2m+1)y_0Y_{m+1}(x_0,y_0)}{-(2m+1)y_0X_{m+1}(x_0,y_0)} = \frac{-Y_{m+1}(x_0,y_0)}{X_{m+1}(x_0,y_0)}\,.$$

Similarly, by the Implicit Function Theorem and ~(\ref{3.10}), there are $\delta > 0$, $\varepsilon > 0$, and a function $x = g(y)$
differentiable on $(y_0-\delta,y_0+\delta)$ such that the set
\begin{equation} \notag \begin{split} C_2 : & = \, \{(x,y): Y_m(x,y) = 0, \enskip y \in (y_0-\delta,y_0+\delta), \enskip x \in (x_0-\varepsilon,x_0+\varepsilon)\} \cr
& = \, \{(g(y),y): \enskip y \in (y_0-\delta,y_0+\delta)\} \cr \end{split}
\end{equation}
is the graph of a differentiable function $x = g(y)$ on $(y_0-\delta,y_0+\delta)$, and the equation for the tangent line to the curve $C_2$ at the point $(x_0,y_0)$ is
$$x-x_0 = \beta(y-y_0)\,,$$
where by ~(\ref{3.11}) and ~(\ref{3.12}) we have
$$\beta  = \left( \frac{\partial{Y_m}}{\partial y} \Big / \frac{\partial{Y_m}}{{\partial x}} \right)(x_0,y_0)
= \frac{-(2m+1)y_0Y_{m+1}(x_0,y_0)}{-(2m+1)y_0X_{m+1}(x_0,y_0)} = \frac{Y_{m+1}(x_0,y_0)}{X_{m+1}(x_0,y_0)}\,.$$

Observe that $\alpha = -\beta$. We conclude that the curves $C_1$ and $C_2$ are orthogonal to each other at $(x_0,y_0)$. As $H_m \subset C_1 \bigcap C_2$,
the point $(x_0,y_0) \in H_m$ is an isolated point of $H_m$.

If ~(\ref{3.8}) holds, then we have
\begin{equation}\label{3.13}\begin{split} \frac{\partial{Y_m}}{\partial y}(x_0,y_0) & = \, \frac{Y_m(x_0,y_0)}{y_0} - (2m+1)y_0Y_{m+1}(x_0,y_0) \cr
& = \, -(2m+1)y_0Y_{m+1}(x_0,y_0) \neq 0 \cr \end{split}
\end{equation}
and
\begin{equation}\label{3.14} \begin{split} \frac{\partial{X_m}}{\partial x}(x_0,y_0) & = \, (2m+1)y_0Y_{m+1}(x_0,y_0) - \frac{(2m)Y_m(x_0,y_0)}{y_0} \cr
& = \, (2m+1)y_0Y_{m+1}(x_0,y_0) \neq 0\,. \end{split}
\end{equation}
Also,
\begin{equation}\label{3.15} \begin{split} \frac{\partial{X_m}}{\partial y}(x_0,y_0) =  \frac{\partial{Y_m}}{\partial x}(x_0,y_0) = -(2m+1)y_0X_{m+1}(x_0,y_0) \end{split}
\end{equation}

By the Implicit Function Theorem and ~(\ref{3.13}), there are $\delta > 0$, $\varepsilon > 0$, and a function $y = f(x)$
differentiable on $(x_0-\delta,x_0+\delta)$ such that the set
\begin{equation} \notag \begin{split} C_3 : & = \,  \{(x,y): Y_m(x,y) = 0, \enskip x \in (x_0-\delta,x_0+\delta), \enskip y \in (y_0-\varepsilon,y_0+\varepsilon)\} \cr
& = \, \{(x,f(x)): \enskip x \in (x_0-\delta,x_0+\delta)\} \cr \end{split}
\end{equation}
is the graph of a differentiable function $y = f(x)$ on $(x_0-\delta,x_0+\delta)$, and the equation for the tangent line to the curve $C_3$ at the point $(x_0,y_0)$ is
$$y-y_0 = \alpha(x-x_0)\,,$$
where
$$ \alpha = \left( \frac{\partial{Y_m}}{\partial x} \Big / \frac{\partial{Y_m}}{{\partial y}} \right)(x_0,y_0)
= \frac{-(2m+1)y_0X_{m+1}(x_0,y_0)}{-(2m+1)y_0Y_{m+1}(x_0,y_0)} = \frac{X_{m+1}(x_0,y_0)}{Y_{m+1}(x_0,y_0)}\,.$$
Similarly, by the Implicit Function Theorem and ~(\ref{3.14}), there are $\delta > 0$, $\varepsilon > 0$, and a function $x = g(y)$ differentiable on
$(y_0-\delta,y_0+\delta)$ such that the set
\begin{equation} \notag \begin{split} C_4 : & = \, \{(x,y): X_m(x,y) = 0, \enskip y \in (y_0-\delta,y_0+\delta), \enskip x \in (x_0-\varepsilon,x_0+\varepsilon)\} \cr
& = \, \{(g(y),y): \enskip y \in (y_0-\delta,y_0+\delta)\} \cr \end{split}
\end{equation}
is the graph of a differentiable function $x = g(y)$ on $(y_0-\delta,y_0+\delta)$, and the equation for the tangent line to the curve $C_4$ at the point $(x_0,y_0)$ is
$$x-x_0 = \beta(y-y_0)$$
where by ~(\ref{3.15}) we have
\begin{equation} \notag \beta = \left( \frac{\partial{X_m}}{\partial y} \Big / \frac{\partial{X_m}}{{\partial x}} \right)(x_0,y_0)
= \frac{-(2m+1)y_0X_{m+1}(x_0,y_0)}{(2m+1)y_0Y_{m+1}(x_0,y_0)} = \frac{-X_{m+1}(x_0,y_0)}{Y_{m+1}(x_0,y_0)}\,.
\end{equation}

Observe that $\alpha = -\beta$. We conclude that the curves $C_3$ and $C_4$ are orthogonal to each other at $(x_0,y_0)$. As $H_m \subset C_3 \bigcap C_4$,
the point $(x_0,y_0) \in H_m$ is an isolated point of $H_m$.
\end{proof}

Now by combining Propositions~\ref{union} and \ref{countable}, we conclude

\begin{prop}\label{Tamas} In the Special Case the set zero set of ${\mathbf F}=(X,Y)$ is countable.
\end{prop}

\begin{cor}\label{XXXtoo} In the Special Case the zero set of ${\mathbf F}=(X,Y)$ in ${\mathbb R}^2$ cannot contain a non-trivial curve.
\end{cor}

But this does not immediately rule out that in some Special Case there are common zeros $(x_m,y_m)$ tending to infinity. So to get the stronger result than 
Corollary~\ref{XXXtoo} that the zero set is actually finite, we have to analyze the situation a little more. The argument is along the lines of the proof of 
Proposition~\ref{squaring}. To prove Proposition~\ref{FINITE} we proceed by showing that the zero set $X=Y=0$ in ${\mathbb R}^2$ is a real {\it algebraic} variety 
and combining this with Gibson~\cite{Gibson} and Corollary~\ref{XXXtoo}. 

\subsection{Proof of Proposition~\ref{FINITE}: In the Special Case the Zero Set of ${\mathbf F}=(X,Y)$ Is Finite}\label{subsec3.4}

\begin{proof}~\ref{FINITE} Again, suppose $0 < x_1 < x_2 < \cdots < x_M$ and the real numbers $a_1, a_2, \ldots, a_M$ are not all zero.
As before, let 
$$X(x,y) = \sum_{j=1}^M{\frac{a_j(x-x_j)}{((x-x_j)^2 + y^2)^{3/2}}}\,, \quad \text{and} 
\quad  Y(x,y) = \sum_{j=1}^M{\frac{a_jy}{((x-x_j)^2 + y^2)^{3/2}}}\,.$$
Let $\Sigma$ be the collection of the $2^M$ functions $\sigma: \{1,2,\ldots,M\} \rightarrow \{-1,1\}$.
Let 
$$X_{\sigma}(x,y) = \sum_{j=1}^M{\frac{\sigma(j)a_j(x-x_j)}{((x-x_j)^2 + y^2)^{3/2}}}\,, \quad \text{and} 
\quad Y_{\sigma}(x,y) = \sum_{j=1}^M{\frac{\sigma(j)a_jy}{((x-x_j)^2 + y^2)^{3/2}}}\,.$$ 
So $(X_{\sigma},Y_{\sigma})$ is also an electric field in the plane with charges on the positive $x$-axis.

Let
$$D_j(x,y) = \prod_{k=1, k \neq j}^M{((x-x_k)^2 + y^2))^{3/2}}\,, \qquad j=1,2,\ldots,M\,,$$
and
$$D(x,y) = \prod_{j=1}^M{((x-x_j)^2 + y^2))^{3/2}}\,.$$
We have
\begin{equation} \notag \begin{split}
X_{\sigma}(x,y)^2 + Y_{\sigma}(x,y)^2 = & \, \Bigg( \sum_{j=1}^M{\frac{\sigma(j)a_j(x-x_j)D_j(x,y)}{D(x,y)}} \Bigg)^2 
+ \Bigg( \sum_{j=1}^M{\frac{\sigma(j)a_jyD_j(x,y)}{D(x,y)}} \Bigg)^2 \cr   
= & \, \frac{1}{D(x,y)^2} \Bigg( \Bigg( \sum_{j=1}^M{\sigma(j)a_j(x-x_j)D_j(x,y)} \Bigg)^2 
+ \Bigg( \sum_{j=1}^M{\sigma(j)a_jyD_j(x,y)} \Bigg)^2 \Bigg)\,. \cr
\end{split} \end{equation} 
Observe that the function 
$$H(x,y) = \prod_{\sigma \in \Sigma}{\Bigg( \Bigg( \sum_{j=1}^M{\sigma(j)a_j(x-x_j)D_j(x,y)} \Bigg)^2 
+ \Bigg( \sum_{j=1}^M{\sigma(j)a_jyD_j(x,y)} \Bigg)^2 \Bigg)}$$
is an even polynomial in each of the variables 
$$\chi_j = D_j(x,y)\,, \qquad j=1,2,\ldots,M\,,$$ 
as it remains the same when $\chi_j$ is replaced by $-\chi_j$. Hence $H(x,y)$ is a polynomial in each of the variables 
$$\chi_j^2 = D_j(x,y)^2, \qquad j=1,2,\ldots,M\,.$$ 
We conclude that 
$$G(x,y) = \prod_{\sigma \in \Sigma}{\Bigg( X_{\sigma}(x,y)^2 + Y_{\sigma}(x,y)^2 \Bigg)} =  \frac{P(x,y)}{D(x,y)^{2^{M+1}}}\,,$$
where $P$ is a polynomial of degree at most $(3M)2^M$.  
We claim that the set 
$$E = \{(x,y): P(x,y) = 0\}$$ 
cannot contain a non-trivial curve. Indeed, if it did then Proposition~\ref{XXX} applied to the finitely many fields $(X_\sigma,Y_\sigma)$ would give a contradiction.
Thus $E = \{(x,y): P(x,y) = 0\}$ is an algebraic variety not containing a non-trivial curve. 
Such an algebraic variety contains only finitely many points.  See Gibson~\cite{Gibson}.   As the zero set of ${\mathbf F}=(X,Y)$ is a 
subset of $E$, it is also finite. 
\end{proof}

\subsection{Upper Bound for the Number of Points in the Zero Set of ${\mathbf F}(X,Y)$}\label{subsec3.5}

With the help of B\'ezout's Theorem we can give an upper bound for the number of points in the
zero set of ${\mathbf F}=(X,Y)$ above.

\noindent To prove the above proposition we need the following classic version of B\'ezout's Theorem.
Two polynomials $f \in {\mathbb R}[x,y]$ and $g \in {\mathbb R}[x,y]$ are coprime if there is no not identically constant polynomial which is 
a factor of both $f$ and $g$. Associated with a polynomial $f \in {\mathbb R}[x,y]$ we define its zero set 
$$Z(f) = \{(x,y): \enskip f(x,y) = 0\}\,.$$

\begin{prop} {\bf (B\'ezout's Theorem)} \label{Classic} Given any two coprime polynomials $f,g \in {\mathbb R}[x,y]$ of degrees $d_1$ and $d_2$, respectively, 
$Z(f) \cap Z(g)$ contains at most $d_1d_2$ points.
\end{prop}

\noindent A consequence of the above B\'ezout's Theorem is the following.

\begin{prop}\label{Bezoutgives} Let $f \in {\mathbb R}[x,y]$ be a polynomial of degree $d$. If $Z(f)$ contains only finitely many points, 
then $Z(f)$ contains at most $d(d-1)$ points.
\end{prop}

Let $f_x$ and $f_y$ be the partial derivatives of $f \in {\mathbb R}[x,y]$ with respect to $x$ and $y$, respectively. To prove Proposition~\ref{Bezoutgives}, we need the following observation.

\begin{prop}\label{Step} Suppose $f  \in {\mathbb R}[x,y]$ is square-free in the sense that $f$ is a product of irreducible polynomial factors none of which is repeated. 
Suppose $P \in {\mathbb R}[x,y]$ is a not identically constant irreducible factor of $f$.
Then either $P$ does not divide $f_x$ or else $P$ does not divide $f_y$.
\end{prop}
\begin{proof} Suppose $P \in {\mathbb R}[x,y]$ is a not identically constant irreducible factor of $f$. Then either $P_x$ is not identically zero 
or else $P_y$ is not identically zero. Without loss of generality we may assume that $P_x$ is not identically zero. Suppose to the contrary that 
$P$ divides $f_x$. We have $f = PQ$ and $f_x = PQ_x + P_xQ$ with some $Q \in {\mathbb R}[x,y]$. Observe that $P$ divides $P_xQ$ and $P$ does not divide $P_x$. 
Since $P \in {\mathbb R}[x,y]$ is irreducible, $P$ divides $Q$. But then $P^2$ divides $f = PQ$. Since $f \in {\mathbb R}[x,y]$ is
square-free this is possible only if $P$ is identically constant, a contradiction.
\end{proof}

\noindent Now we are ready to prove Proposition~\ref{Bezoutgives}.

\begin{proof} Without loss of generality we may assume that $f$ is square-free, that is, $f$ is a product of irreducible polynomial factors none of which is repeated,
as keeping only one of each repeated factors in a factorization of $f$ the set $Z(f)$ remains the same.
Observe that if $Z(f)$ contains only finitely many points, then each of these points are singular, that is, $Z(f) \subset Z(f_x)$ and $Z(f) \subset Z(f_y)$,
otherwise the Implicit Function Theorem would imply that $Z(f)$ contains a curve (graph of a function) with infinitely many points. Let 
$f = P_1P_2 \cdots P_m$, where the polynomials $P_j \in {\mathbb R}[x,y]$ are not identically constant and irreducible. Denote the degree of 
$P_j$ by $d_j$ and denote the degree of $P$ by $d$. We have $\sum_{j=1}^m{d_j} = d$. By Proposition~\ref{Step} for every fixed $j$ either 
$P_j$ and $f_x$ are coprime or $P_j$ and $f_y$ are coprime. If $P_j$ and $f_x$ are coprime, then B\'ezout's Theorem implies
$$Z(P_j) = Z(P_j) \cap Z(f) \subset Z(P_j) \cap Z(f_x)$$ has at most $d_j(d-1)$ points. If $P_j$ and $f_y$ are coprime, then B\'ezout's Theorem implies
$$Z(P_j) = Z(P_j) \cap Z(f) \subset Z(P_j) \cap Z(f_y)$$ has at most $d_j(d-1)$ points. Summing for $j = 1,2,\ldots,m$ we obtain that
$Z(f)$ contains at most $d(d-1)$ points.
\end{proof}

\noindent Now we are ready to prove Proposition~\ref{UPPERBOUND}.

\begin{proof} Observe that the degree of $P$ in the proof of Proposition~\ref{FINITE} is at most $d = (3M)2^M$. Hence, by Proposition~\ref{Bezoutgives} we can deduce that 
the zero set $Z(P) = \{(x,y): P(x,y) = 0\}$ contains at most $d^2 = 9M^24^M$ points. As the zero set of ${\mathbf F}=(X,Y)$ is a subset of $Z(P)$, 
it contains at most $d^2 = 9M^24^M$ points as well.
\end{proof}

Observe that no general fact about analytic varieties has been used. All that is needed is B\'ezout's Theorem. See Gibson~\cite{Gibson} for B\'ezout's Theorem.

\section{Special Case: Asymptotic Directions of Component Zero Sets}\label{Asymptotes}

The analysis in this section followed after computing many examples of the Special Case using Mathematica.  It is clear that quite often there are zeros of such electric fields not on the axis of the charges.  Also, in most cases (but possibly not all), when there are zeros of $X$ and $Y$ going to infinity, then the zero sets of $X$ and $Y$ separately lie on curves (perhaps always are actually curves) that are asymptotic to lines going to infinity, in directions that are interlacing.  

So we sought methods to prove that the asymptotics at infinity are indeed as the computed examples suggested.  This ought to be not so difficult, but we have had to resort to the somewhat complicated arguments in this section.  Moreover, because there is a need for series expansions, it was necessary to divide the cases into types.  
In any case, this analysis leads to interesting sequences of polynomials whose zeros give the possible asymptotic directions at infinity of the zero sets of $X$ and $Y$.   

\begin{defn}\label{AD1} 
Let $\alpha \in (-\infty,\infty)$. The line $y = \alpha x$ is called an {\em asymptotic direction} for a set $A \subset {\mathbb R}^2$ 
if there are $(p_m,q_m) \in A$ such that $|p_m| \rightarrow \infty$ and 
$$\lim_{|m \rightarrow \infty}{\frac{q_m}{p_m}} = \alpha\,.$$
\end{defn}

\begin{defn}\label{ADy}
The $y$-axis is called an asymptotic direction for a set $A \subset {\mathbb R}^2$ 
if there are $(p_m,q_m) \in A$ such that $|q_m| \rightarrow \infty$ and 
$$\lim_{m \rightarrow \infty}{\frac{p_m}{q_m}} = 0\,.$$  
\end{defn}

For the sake of brevity we will use the notation $\beta = 1/\alpha$. Note that the case $\beta = 0$ corresponds to the asymptotic direction given by the $y$-axis.
The goal in this section is simple: to show that the set of possible asymptotic directions for the zero set $X=0$ and the set of the possible asymptotic directions 
for the zero set $Y=0$ are distinct. We succeed in showing this with the exception of the asymptotic directions $y = \pm x$.  

\begin{defn}\label{Types} We call the domain $\{(x,y):|x| < |y|\}$ Type I Domain. We call $\{(x,y):|y| < |x|\}$ Type II Domain. 
\end{defn}
In fact, for technical reasons, for a fixed $\delta \in (0,1/2)$, we will consider Type I Domains
$\{(x,y):|x| < (1-\delta)|y|\}$ and Type II Domains $\{(x,y):|y| < (1-\delta)|x|\}$ nd we let $\delta > 0$ tend to $0$. 

\noindent To find all possible asymptotic directions for the zero sets $Y=0$ and $X=0$ we proceed differently in Type I and Type II Domains, but the cases of 
the zero sets $Y=0$ and $X=0$ are quite similar to each other in both domains.  
We think that our results in this section have their own intrinsic interest, even though the boundedness of the zero set $X=Y=0$ that obviously follows from 
Proposition~\ref{FINITE}, is not completely recaptured in this section as we cannot show that the lines $|y|=|x|$ cannot be a common asymptotic direction 
to both of the zero sets $Y=0$ and $X=0$. 

By the Binomial Theorem the Taylor series expansion
$$(1+z)^{-3/2} = \sum_{n=0}^\infty{\binom {-3/2}{n}{z^n}}\,, \qquad z \in (-1,1)\,,$$
converges uniformly on the interval $[-1+\delta,1-\delta]$ for every $\delta > 0$.
Let $\mu_u = (-1)^u \sum_{j=1}^M{a_jx_j^u}$.

Recall Definitions~\ref{SC} and \ref{L}. 
This is the choice of $L$ throughout Sections~\ref{SecYTypeI}, \ref{SecXTypeI}, \ref{SecYTypeII}, and \ref{SecXTypeII}. Typically, as $L$ gets larger, 
the complexity to detect possible asymptotic directions for the zero sets $X=0$ and $Y=0$ increases.

\subsection{Type I Domain, Possible Asymptotic Directions for the Zero Set $Y=0$}\label{SecYTypeI}

In this section we present the proof of Proposition~\ref{YTypeI}.

\begin{proof}
We have
\begin{equation} \notag
\begin{split}
0 & = \,  q_m^{L+1}|q_m|Y(p_m,q_m) =  q_m^{L}\sum_{j=1}^M{a_j \left(1 + \left( \frac{p_m-x_j}{q_m} \right)^2 \right)^{-3/2}} \cr
& = \, q_m^{L}\sum_{j=1}^M{a_j \sum_{n=0}^\infty{\binom{-3/2}{n} \left( \frac{p_m-x_j}{q_m} \right)^{2n}}} \cr
\end{split}
\end{equation}
for all sufficiently large $m$. Note that we have adjusted with the factor $q_m^{L+1}|q_m|$. We also have 
\begin{equation} \label{F1} 
\begin{split}
\sum_{j=1}^M & {a_j(p_m-x_j)^{2n}} =\\ 
\quad & \sum_{j=1}^M{a_j \sum_{u=0}^{2n}{\binom{2n}{u}(-1)^ux_j^up_m^{2n-u}}} 
= \sum_{u=0}^{2n}{\sum_{j=1}^M{\binom{2n}{u}(-1)^ua_jx_j^up_m^{2n-u}}} \cr 
& = \, \sum_{u=L}^{2n}{\sum_{j=1}^M{\binom{2n}{u}(-1)^ux_j^u p_m^{2n-u}}} = \sum_{u=L}^{2n}{\binom {2n}{u} \mu_u p_m^{2n-u}}\,, \cr 
\end{split}
\end{equation}
and hence by (\ref{F1})
\begin{equation} \notag \begin{split}
0 = & q_m^{L}\sum_{j=1}^M{a_j \sum_{n=0}^\infty{\binom{-3/2}{n} \left( \frac{p_m-x_j}{q_m} \right)^{2n}}} 
= q_m^{L}\sum_{n=0}^\infty{ \binom{-3/2}{n} \sum_{j=1}^M{a_j \left( \frac{p_m-x_j}{q_m} \right)^{2n}}} \cr
= &  \sum_{n=0}^\infty{\binom {-3/2}{n} \sum_{u=L}^{2n}{\binom{2n}{u} \mu_u \left( \frac{p_m}{q_m} \right)^{2n-u} q_m^{L-u}}} \cr
\end{split}\end{equation}
for all sufficiently large $m$. Separating the term for which $u=L$, we obtain  
\begin{equation}\label{F2} \begin{split}
0 = & \sum_{n=0}^\infty{\binom {-3/2}{n} {\binom{2n}{L} \mu_L \left( \frac{p_m}{q_m} \right)^{2n-L}}} 
+ \sum_{n=0}^\infty{\binom {-3/2}{n} \sum_{u=L+1}^{2n}{\binom{2n}{u} \mu_u \left( \frac{p_m}{q_m} \right)^{2n-u} q_m^{L-u}}} \cr 
= & \, \frac{\mu_L}{L!} \frac{d^L}{dz^L} \left( \left( 1 + z^2 \right)^{-3/2} \right) \Big|_{p_m/q_m}
+ \sum_{u=L+1}^\infty{\frac{\mu_u}{u!} \frac{d^u}{dz^u} \left( \left( 1 + z^2 \right)^{-3/2} \right) \Big|_{p_m/q_m}} q_m^{L-u} \cr
\end{split}\end{equation}
for all sufficiently large $m$. Here we interchanged the order of summations, which is legal because one can easily check that under our conditions
the double sum converges absolutely. See Section~\ref{interchange}. 
Observe that
\begin{equation}\label{F3}
|\mu_u| = \left| \sum_{j=1}^M{a_jx_j^u} \right| \leq \left( \sum_{j=1}^M{|a_j|} \right) x_M^u = A x_M^u\,.
\end{equation}
The Cauchy Integral Formula and $|p_m/q_m| < 1-\delta$ imply that
\begin{equation}\label{F4}
\Bigg| \frac{1}{u!} \frac{d^u}{dz^u}\left( 1 + z^2 \right)^{-3/2} 
 \Big|_{p_m/q_m} \Bigg| 
\leq \frac{\delta}{2} \, \max_{|z| = 1 - \delta/2} {\Bigg| \left( 1 + z^2 \right)^{-3/2} \Bigg|} \left( \frac{\delta}{2} \right)^{-(u+1)} \leq \left( \frac{2}{\delta} \right)^{u+2}\,.
\end{equation}
Observe also that
\begin{equation}\label{F5}
|q_m|^{L-u} = |q_m|^{L-u+1/2} |q_m|^{-1/2} \leq |q_m|^{-u/(2L+2)} |q_m|^{-1/2}\, \quad u \geq L+1\,, \enskip |q_m| > 1\,. 
\end{equation}
Combining  (\ref{F3}), (\ref{F4}), and (\ref{F5}),  we obtain
\begin{equation} \label{F6}
\begin{split}
& \, \Bigg| \sum_{u=L+1}^\infty{\frac{\mu_u}{u!} \frac{d^u}{dz^u} \left( \left( 1 + z^2 \right)^{-3/2} \right) \Big|_{p_m/q_m}} q_m^{L-u} \Bigg| \cr 
\leq & \, \sum_{u=L+1}^\infty{ |\mu_u| \Bigg| \frac{1}{u!} \frac{d^u}{dz^u} \left( \left( 1 + z^2 \right)^{-3/2} \right) \Big|_{p_m/q_m} \Bigg| |q_m|^{L-u}} \cr 
\leq & \, \sum_{u=L+1}^\infty{A x_M^u \left( \frac{2}{\delta} \right)^{u+2}  |q_m|^{L-u}} 
\leq |q_m|^{-1/2} \sum_{u=L+1}^\infty{A x_M^u \left( \frac{2}{\delta} \right)^{u+2} |q_m|^{-u/(2L+2)}} \cr
\leq & \, |q_m|^{-1/2} A\left( \frac{2}{\delta} \right)^2 \sum_{u=L+1}^\infty{\left( \frac{2x_M}{\delta |q_m|^{1/(2L+2)}} \right)^u}\,. \cr 
\end{split}
\end{equation}
As $|q_m| \rightarrow \infty$, (\ref{F6}) implies that
\begin{equation} \notag
\lim_{m \rightarrow \infty}{ \left( \sum_{u=L+1}^\infty{\frac{\mu_u}{u!} \frac{d^u}{dz^u} \left( \left( 1 + z^2 \right)^{-3/2} \right) \Big|_{p_m/q_m}} q_m^{L-u} \right)} = 0\,.
\end{equation}
Now we are ready to take the limit in (\ref{F2}) when $|q_m| \rightarrow \infty$. We obtain
$$\frac{\mu_L}{L!} \frac{d^L}{d z^L} \left( \frac{1}{(1+z^2)^{3/2}} \right)\bigg |_\beta  + 0 = 0\,.$$
\end{proof}

\subsection{Type I Domain, Possible Asymptotic Directions for the Zero Set $X=0$}\label{SecXTypeI}

In this section we present the proof of Proposition~\ref{XTypeI}, which is very similar to the proof of Proposition~\ref{YTypeI}.

\begin{proof}
We have
\begin{equation} \notag
\begin{split}
0 & =  \,  q_m^{L+1}|q_m| X(p_m,q_m) =  q_m^{L-1}\sum_{j=1}^M{a_j(p_m-x_j)\left(1 + \left( \frac{p_m-x_j}{q_m} \right)^2 \right)^{-3/2}} \cr
& = \, q_m^{L-1} \sum_{j=1}^M{a_j \sum_{n=0}^\infty{\binom{-3/2}{n} \frac{(p_m-x_j)^{2n+1}}{q_m^{2n}}}} \cr
\end{split}
\end{equation}
for all sufficiently large $m$. Note that we have adjusted with the factor $q_m^{L+1}|q_m|$. Replacing $2n$ by $2n+1$ in (\ref{F2}), we also have
\begin{equation} \notag
\sum_{j=1}^M{a_j(p_m-x_j)^{2n+1}} = \sum_{u=L}^{2n+1}{\binom {2n+1}{u} \mu_u p_m^{2n+1-u}}\,,
\end{equation}
and hence
\begin{equation} \notag \begin{split}
0 & = \, q_m^{L-1} \sum_{j=1}^M{a_j \sum_{n=0}^\infty{\binom{-3/2}{n} \frac{(p_m-x_j)^{2n+1}}{q_m^{2n}}}}
= q_m^{L-1}\sum_{n=0}^\infty{ \binom{-3/2}{n} \sum_{j=1}^M{ \frac{(p_m-x_j)^{2n+1}}{q_m^{2n}}}} \cr
& = \, \sum_{n=0}^\infty{\binom {-3/2}{n} \sum_{u=L}^{2n+1}{\binom{2n+1}{u} \mu_u \left( \frac{p_m}{q_m} \right)^{2n+1-u} q_m^{L-u}}}\,. \cr
\end{split}\end{equation}
for all sufficiently large $m$. Separating the term for which $u=L$, we obtain
\begin{equation} \label{F7}
\begin{split}
& \sum_{n=0}^\infty{\binom {-3/2}{n} {\binom{2n+1}{L} \mu_L \left( \frac{p_m}{q_m} \right)^{2n+1-L}}}
+ \sum_{n=0}^\infty{\binom {-3/2}{n} \sum_{u=L+1}^{2n+1}{\binom{2n+1}{u} \mu_u \left( \frac{p_m}{q_m} \right)^{2n+1-u} q_m^{L-u}}} \cr
= & \, \frac{\mu_L}{L!} \frac{d^L}{dz^L} \left( z \left( 1 + z^2 \right)^{-3/2} \right) \Big|_{z = p_m/q_m}
+ \sum_{u=L+1}^\infty{\frac{\mu_u}{u!} \frac{d^u}{dz^u} \left( z \left( 1 + z^2 \right)^{-3/2} \right) \Big|_{z = p_m/q_m}} q_m^{L-u} \cr
\end{split}
\end{equation}
vanishes for all sufficiently large $m$. Here we interchanged the order of summations, which is legitimate as one can easily check that under
our conditions the double sum converges absolutely. See Section~\ref{interchange}.
The Cauchy Integral Formula and $|p_m/q_m| \leq 1-\delta$ imply that
\begin{equation} \label{F8}
\begin{split}
\Bigg| \frac{1}{u!} \frac{d^u}{dz^u} \left( z\left( 1 + z^2 \right)^{-3/2} \right) \Big|_{z = p_m/q_m} \Bigg|
\leq \frac{\delta}{2} \, \max_{|z| = 1 - \delta/2} {\Bigg| \left( 1 + z^2 \right)^{-3/2} \Bigg|} \left( \frac{\delta}{2} \right)^{-(u+1)} \leq \left( \frac{2}{\delta} \right)^{u+2}\,.
\end{split}
\end{equation}
Combining (\ref{F3}), (\ref{F7}), and (\ref{F8}), we obtain
\begin{equation} \label{F9}
\begin{split}
& \, \Bigg| \sum_{u=L+1}^\infty{\frac{(-1)^u \mu_u}{u!} \frac{d^u}{dz^u} \left( z\left( 1 + z^2 \right)^{-3/2} \right) \Big|_{z = p_m/q_m}} q_m^{L-u} \Bigg| \cr
\leq & \, \sum_{u=L+1}^\infty{ |\mu_u| \Bigg| \frac{1}{u!} \frac{d^u}{dz^u} \left(z \left( 1 + z^2 \right)^{-3/2} \right) \Big|_{z = p_m/q_m} \Bigg| |q_m|^{L-u}} \cr
\leq & \, \sum_{u=L+1}^\infty{A x_M^u \left( \frac{2}{\delta} \right)^{u+2}  |q_m|^{L-u}}
\leq |q_m|^{-1/2} \sum_{u=L+1}^\infty{A x_M^u \left( \frac{2}{\delta} \right)^{u+2} |q_m|^{-u/(2L+2)}} \cr
\leq & \, |q_m|^{-1/2} A\left( \frac{2}{\delta} \right)^2 \sum_{u=L+1}^\infty{\left( \frac{2x_M}{\delta|q_m|^{1/(2L+2)}} \right)^u}\,. \cr
\end{split}
\end{equation}
As $m \rightarrow \infty$, (\ref{F9}) implies that
\begin{equation} \notag
\lim_{m \rightarrow \infty}{ \left( \sum_{u=L+1}^\infty{\frac{\mu_u}{u!} \frac{d^u}{dz^u} \left( z\left( 1 + z^2 \right)^{-3/2} \right) \Big|_{z = p_m/q_m}}q_m^{L-u} \right)} = 0\,.
\end{equation}
Now we are ready to take the limit in \ref{F7} when $m \rightarrow \infty$. We obtain
$$\frac{\mu_L}{L!} \frac{d^L}{d z^L} \left( \frac{\beta}{(1+z^2)^{3/2}} \right) \bigg |_\beta + 0 = 0\,.$$
\end{proof}                                                                                                                                                                                                                                                                                                                                                                                                             \subsection{Type II Domain, Possible Asymptotic Directions for the Zero Set $Y=0$}\label{SecYTypeII}

In this section we present the proof of Proposition~\ref{YTypeII}, which is similar to the proof of Proposition~\ref{YTypeI}, but there are some important technical differences. 

\begin{proof}
The Binomial Theorem gives for $|x| > x_M$,
\begin{equation} \notag 
\frac{1}{(x-x_j)^{2n+3}} = \, \frac{1}{x^{2n+3}}\left( 1 - \frac{x_j}{x} \right)^{-(2n+3)} 
= \frac{1}{x^{2n+3}} \sum_{u=0}^\infty{\binom{-(2n+3)}{u} \left( \frac{-x_j}{x} \right)^u}\,, 
\end{equation}
for each $j=1,2,\ldots,M$, and hence
\begin{equation}\label{F10} 
\begin{split} 
\sum_{j=1}^M{\frac{a_j}{(x-x_j)^{2n+3}}} & = \, \sum_{j=1}^M \frac{a_j}{x^{2n+3}} \sum_{u=0}^\infty{\binom{-(2n+3)}{u} \left( \frac{-x_j}{x} \right)^u} \cr
& = \, \sum_{j=1}^M \frac{a_j}{x^{2n+3}} \left( \sum_{u=0}^L{\binom{-(2n+3)}{u} \left( \frac{-x_j}{x} \right)^u} + R_{L,n}\left( \frac{-x_j}{x} \right) \right) \cr
& = \, \mu_L \binom{-(2n+3)}{L}x^{-(2n+3)-L} + A_{L,n}(x)\,, \qquad |x| > x_M\,, \cr 
\end{split} 
\end{equation}
with
$$A_{L,n}(x) = x^{-(2n+3)}\sum_{j=1}^M{a_j R_{L,n}\left( \frac{-x_j}{x} \right)}\,,$$
where $R_{L,n}(x)$ is the $L$th remainder term in the Taylor series expansion of the function 
$$f(t) := (1+t)^{-(2n+3)}$$ 
centered at $0$, that is, 
$$(1+t)^{-(2n+3)} = T_{L,n}(t) + R_{L,n}(t)\,, \qquad |t| < 1\,,$$
where $T_{L,n}(t)$ it the $L$th Taylor polynomial centered at $0$ associated with the function $f$. 
Estimating by using the Cauchy form of the remainder term $R_{L,n}(t)$ in the Taylor series expansion of the function $f$, we obtain
$$|R_{L,n}(t)| \leq \, (L+1) \left| \binom{-(2n+3)}{L+1} \right| (1-|t|)^{-(2n+4)}|t|^{L+1}\,, \qquad |t| < 1\,,$$ 
and hence with $A = \sum_{j=1}^M{|a_j|}$, we have for $|x| > x_M$,
\begin{equation}\label{F11}
\begin{split} 
|A_{L,n}(x)| & = \, \left| x^{-(2n+3)}\sum_{j=1}^M{a_jR_{L,n}\left( \frac{-x_j}{x} \right)} \right| \cr 
& \leq \, A|x|^{-(2n+3)}(L+1)\left| \binom{-(2n+3)}{L+1} \right| \left(1 - \left| \frac{x_M}{x} \right| \right)^{-(2n+4)} \left| \frac{x_M}{x} \right|^{L+1} \, \cr
& \leq \, A(L+1)x_M^{L+1} \left| \binom{-(2n+3)}{L+1} \right| \left( |x| - x_M \right)^{-(2n+4)} |x|^{-L}\, . \cr
\end{split} 
\end{equation}
Some simple algebra shows that for $|x| > x_M, y \neq 0$,
$$\frac{Y(x,y)}{y} = \sum_{j=1}^M{\frac{a_j}{((x-x_j)^2 + y^2)^{3/2}}} = 
\sum_{j=1}^M{\frac{a_j}{|x-x_j|^3} \left( 1 +  \left( \frac{y}{x-x_j} \right)^2 \right)^{-3/2}}\, .$$
Observe that $\displaystyle{\left| \frac{q_m}{p_m-x_j} \right| \leq 1-\delta}$ for all sufficiently large $m$, and hence the Binomial Theorem gives, for $q_m\not= 0$,
\begin{equation} \notag \begin{split}
0 & = \, \pm \frac{p_m^{L+3}}{q_m} Y(p_m,q_m) = p_m^{L+3} \sum_{j=1}^M{\frac{a_j}{(p_m-x_j)^3} \left( 1 +  \left( \frac{q_m}{p_m-x_j} \right)^2 \right)^{-3/2}} \cr 
& = \, p_m^{L+3} \sum_{j=1}^M{\frac{a_j}{(p_m-x_j)^3} \sum_{n=0}^\infty{\binom{-3/2}{n} \left( \frac{q_m}{p_m-x_j} \right)^{2n}}} \cr
\end{split} \end{equation}
for all sufficiently large $m$. Note that we have adjusted with the factor $p_m^{L+3}$. Combining this with (\ref{F10}), we obtain 
\begin{equation}\label{F12} 
\begin{split} 0 = & \, p_m^{L+3} \sum_{j=1}^M{\frac{a_j}{(p_m-x_j)^3} \sum_{n=0}^\infty{\binom {-3/2}{n} \left( \frac{q_m}{p_m-x_j} \right)^{2n}}} \cr 
= & \, p_m^{L+3} \sum_{n=0}^\infty{\sum_{j=1}^M{\binom{-3/2}{n} \frac{a_j q_m^{2n}}{(p_m-x_j)^{2n+3}}}} \cr
= & \, \sum_{n=0}^\infty{\binom {-3/2}{n} q_m^{2n} \left( \mu_L \binom{-(2n+3)}{L}p_m^{-2n} + p_m^{L+3}A_{L,n}(p_m) \right)} \cr 
= & \, \sum_{n=0}^\infty{\binom {-3/2}{n} \left( \mu_L \binom{-(2n+3)}{L} \left( \frac{q_m}{p_m} \right)^{2n} + q_m^{2n}p_m^{L+3}A_{L,n}(p_m) \right)}
\end{split} \end{equation}
for all sufficiently large $m$. With
$$ F(m,M) = \left( \frac{|p_m|}{|p_m|-x_M} \right)^3 \frac{1}{|p_m|-x_M}  , $$
it follows from (\ref{F11}) that
\begin{equation} \notag \begin{split}  
&  \, \left| \sum_{n=0}^\infty{ \binom {-3/2}{n} q_m^{2n}p_m^{L+3}A_{L,n}(p_m)} \right| \cr   
\leq & \, A(L+1)x_M^{L+1} \sum_{n=0}^\infty{ \left| \binom {-3/2}{n} \binom{-(2n+3)}{L+1} \right| \left( \frac{|q_m|}{|p_m|-x_M} \right)^{2n}  F(m,M)} \cr 
\leq & \, A(L+1)x_M^{L+1} \sum_{n=0}^\infty{ \left| \binom {-3/2}{n} \binom{-(2n+3)}{L+1} \right| \left( |\alpha| + \delta \right)^{2n} F(m,M)}
\end{split} \end{equation}
for all sufficiently large $m$, where the right-hand side converges to
$$A(L+1)x_M^{L+1} \left( \sum_{n=0}^\infty{ \left| \binom {-3/2}{n} \binom{-(2n+3)}{L+1} \right| \left( |\alpha| + \delta \right)^{2n}}  \right) 
\left( \lim_{m  \rightarrow \infty}{ \frac{|p_m|^3}{(|p_m|-x_M)^4}} \right).$$
Since this limit is zero, taking the limit in (\ref{F12}) as $m \rightarrow \infty$ gives
$$\sum_{n=0}^\infty{\binom{-3/2}{n} \mu_L \binom{-(2n+3)}{L} \alpha^{2n}} = 0\,.$$
That is, 
$$\sum_{n=0}^\infty{\binom{-3/2}{n} \mu_L \binom{2n+L+2}{L} \alpha^{2n}} = 0\,.$$
As $|\alpha| \leq 1-\delta$, we conclude that $\alpha \neq 0$ and multiplying by $\alpha^2$ gives
$$\frac{\mu_L}{L!} \frac{d^L}{dz^L} \left( \frac{z^{L+2}}{(1+z^2)^{3/2}} \right)\bigg |_\alpha = 0\,.$$
\end{proof}

\subsection{Type II Domain, Possible Asymptotic Directions for the Zero Set $X=0$}\label{SecXTypeII}

In this section we present the proof of Proposition~\ref{XTypeII}, which is similar to the proof of Proposition~\ref{YTypeII}.

\begin{proof}
Replacing $2n+3$ by $2n+2$ in (\ref{F10}), we obtain
\begin{equation} \label{F13}
\sum_{j=1}^M{\frac{a_j}{(x-x_j)^{2n+2}}} = \mu_L \binom{-(2n+2)}{L}x^{-(2n+3)-L} + A_{L,n}(x)\,, \qquad |x| > |x_M|\,,
\end{equation}
with
$$A_{L,n}(x) := x^{-(2n+2)}\sum_{j=1}^M{a_j R_{L,n}\left( \frac{-x_j}{x} \right)}\,,$$
where $R_{L,n}(x)$ is the $L$th remainder term in the Taylor series expansion of the function 
$$f(t) := (1+t)^{-(2n+2)}$$ 
centered at $0$, that is,
$$(1+t)^{-(2n+2)} = T_{L,n}(t) + R_{L,n}(t)\,, \qquad |t| < 1\,,$$
where $T_{L,n}(t)$ it the $L$th Taylor polynomial centered at $0$ associated with the function $f$.
Estimating by using the Cauchy form of the remainder term $R_{L,n}(t)$ in the Taylor series expansion of the function $f$, we obtain
$$|R_{L,n}(t)| \leq \, (L+1) \left| \binom{-(2n+2)}{L+1} \right| (1-|t|)^{-(2n+3)}|t|^{L+1}\,, \qquad |t| < 1\,,$$
and hence, replacing $2n+3$ with $2n+2$ in (\ref{F11}), we have
\begin{equation}\label{F14}
|A_{L,n}(x)| \leq AM(L+1)x_M^{L+1} \left| \binom{-(2n+2)}{L+1} \right| \left( |x| - x_M \right)^{-(2n+3)} |x|^{-L}\,, \qquad |x| > x_M\,.
\end{equation}
A simple algebra shows that
$$X(x,y) = \sum_{j=1}^M{\frac{a_j(x-x_j)}{((x-x_j)^2 + y^2)^{3/2}}} =
\sum_{j=1}^M{\frac{a_j(x-x_j)}{|x-x_j|^3} \left( 1 +  \left( \frac{y}{x-x_j} \right)^2 \right)^{-3/2}}\,, \quad |x| > x_M, \enskip y \neq 0\,.$$
Observe that
$$\left| \frac{q_m}{p_m-x_j} \right| \leq 1-\delta$$
for all sufficiently large $m$, and hence the Binomial Theorem gives
\begin{equation} \notag \begin{split}
0 & = \, \pm p_m^{L+2}X(p_m,q_m) = p_m^{L+2} \sum_{j=1}^M{\frac{a_j(p_m-x_j)}{(p_m-x_j)^3} \left( 1 +  \left( \frac{q_m}{p_m-x_j} \right)^2 \right)^{-3/2}} \cr
& = \, p_m^{L+2} \sum_{j=1}^M{\frac{a_j}{(p_m-x_j)^2} \sum_{n=0}^\infty{\binom{-3/2}{n} \left( \frac{q_m}{p_m-x_j} \right)^{2n}}} \cr
\end{split} \end{equation}
for all sufficiently large $m$. Note that we have adjusted with the factor $p_m^{L+2}$. Combining this with (\ref{F13}), similarly to (\ref{F12}), we get
\begin{equation} \label{F15}
\begin{split} 0 & = \, p_m^{L+2} \sum_{j=1}^M{\frac{a_j}{(p_m-x_j)^2} \sum_{n=0}^\infty{\binom {-3/2}{n} \left( \frac{q_m}{p_m-x_j} \right)^{2n}}} \cr
& = \, \sum_{n=0}^\infty{\binom {-3/2}{n} \left( \mu_L \binom{-(2n+2)}{L} \left( \frac{q_m}{p_m} \right)^{2n} + q_m^{2n}p_m^{L+2}A_{L,n}(p_m) \right)}
\end{split} \end{equation}
for all sufficiently large $m$. It follows from (\ref{F14}) that
\begin{equation} \notag \begin{split}
&  \left| \sum_{n=0}^\infty{\binom {-3/2}{n} q_m^{2n}p_m^{L+2}A_{L,n}(p_m)} \right| \cr
\leq & AM(L+1)x_M^{L+1} \sum_{n=0}^\infty{ \left| \binom {-3/2}{n} \binom{-(2n+2)}{L+1} \right| \left( \frac{|q_m|}{|p_m|-x_M} \right)^{2n}
\left( \frac{|p_m|}{|p_m|-x_M} \right)^2 \frac{1}{|p_m| - x_M}} \cr
\leq & AM(L+1)x_M^{L+1} \sum_{n=0}^\infty{ \left| \binom {-3/2}{n} \binom{-(2n+2)}{L+1} \right| \left( |\alpha| + \delta \right)^{2n}}
\left( \frac{|p_m|}{|p_m|-x_M} \right)^2 \frac{1}{|p_m| -x_M}\cr
\end{split} \end{equation}
for all sufficiently large $m$, where the right-hand side converges to
$$AM(L+1)x_M^{L+1} \left( \sum_{n=0}^\infty{ \left| \binom {-3/2}{n} \binom{-(2n+2)}{L+1} \right| \left( |\alpha| + \delta \right)^{2n}} \right)
\left( \lim_{m \rightarrow \infty}{\frac{|p_m|^2}{(|p_m|-x_M)^3}} \right) = 0\,.$$
Therefore taking the limit in (\ref{F15}) as $m \rightarrow \infty$ gives
$$\sum_{n=0}^\infty{\binom{-3/2}{n} \mu_L \binom{-(2n+2)}{L} \alpha^{2n}} = 0\,.$$
That is,
$$\sum_{n=0}^\infty{\binom{-3/2}{n} \mu_L \binom{2n+L+1}{L} \alpha^{2n}} = 0\,.$$
As $|\alpha| \leq 1-2\delta$ we conclude that $\alpha \neq 0$ and
$$\frac{\mu_L}{L!} \frac{d^L}{d\alpha^L} \left( \frac{\alpha^{L+1}}{(1+\alpha^2)^{3/2}} \right) = 0\,.$$
\end{proof}

\subsection{Interchanging the Order of Summation}\label{interchange}

In the analysis in Section~\ref{SecYTypeI} to find the possible asymptotes for the zero set $Y=0$ in the Type I Domain, we interchanged the order of summations in
(\ref{F2}), that is,
\begin{equation} \notag 
\sum_{n=0}^\infty{\binom{-3/2}{n}{\sum_{u=L}^{2n}{\binom{2n}{u} \mu_u \left( \frac{p_m}{q_m} \right)^{2n-u} q_m^{L-u}}}}  
= \sum_{u=L}^\infty{\sum_{n=0}^\infty{\binom{2n}{u} \mu_u \left( \frac{p_m}{q_m} \right)^{2n-u} q_m^{L-u}}} 
\end{equation}
for all sufficiently large $m$, or equivalently
\begin{equation} \notag 
\sum_{n=0}^\infty{\binom{-3/2}{n}{\sum_{u=L}^{2n}{\binom{2n}{u} \mu_u \left( \frac{p_m}{q_m} \right)^{2n-u} q_m^{-u}}}}  
= \sum_{u=L}^\infty{\sum_{n=0}^\infty{\binom{2n}{u} \mu_u \left( \frac{p_m}{q_m} \right)^{2n-u} q_m^{-u}}} 
\end{equation}
for all sufficiently large $m$. To see that this is legitimate, recall that $\mu_u = (-1)^u\sum_{j=1}^u{a_jx_j^u}$ and 
$\mu_u = 0$ for each $u=0,1,\ldots,L-1$, and hence we have
\begin{equation} \notag 
\begin{split}
&\, \sum_{n=0}^\infty{ \sum_{u=L}^{2n}{\left| \binom {-3/2}{n} \binom{2n}{u} \mu_u \left( \frac{p_m}{q_m} \right)^{2n-u} q_m^{-u} \right|}}  \cr
&= \sum_{n=0}^\infty{ \sum_{u=0}^{2n}{\left| \binom {-3/2}{n} \binom{2n}{u} \mu_u \left( \frac{p_m}{q_m} \right)^{2n-u} q_m^{-u} \right|}} \cr
&= \,  \sum_{n=0}^\infty{\sum_{u=0}^{2n}{\left| \binom {-3/2}{n} \binom{2n}{u} \sum_{j=1}^M{a_j(-x_j)^u} \left( \frac{p_m}{q_m} \right)^{2n-u} q_m^{-u} \right|}}\cr
&\leq \sum_{j=1}^M{\sum_{n=0}^\infty{|a_j| \left| \binom {-3/2}{n} \right| \left( \frac{|p_m| + x_j}{|q_m|} \right)^{2n}}} \cr
&\leq \,  \sum_{j=1}^M{\sum_{n=0}^\infty{|a_j| \left| \binom {-3/2}{n} \right| \left( |\beta| + \delta/2 \right)^{2n}}} < \infty \cr
\end{split}
\end{equation}
for all sufficiently large $m$, as $|\beta| \leq 1-\delta$ and the power series $\sum_{n=0}^\infty{\binom {-3/2}{n} z^{2n}}$ converges absolutely 
for all $z \in (-1,1)$. 

A similar argument shows that interchanging the order of summation in our work in Section~\ref{SecXTypeI} to find the possible asymptotes for the zero set $X=0$ in 
the Type I Domain is legitimate.

\subsection{Interlacing properties of zeros of polynomials}\label{lace}

We now address how we can conclude that the asymptotic directions of the zero sets of $Y$ and $X$ are distinct.  Proving this would show that there are no zeros
for ${\mathbf F}=(X,Y)$ outside some large disk, in the Special Case.  From the analysis in Section~\ref{Curves}, we already know this because the zero set is finite.
But our argument may have independent value because it gives interlacing properties for a class of functions that has not previously been considered. In addition,
this interlacing is needed for applications of the General Case that we will address later.  Moreover, this interlacing sets up a context in which some of the
challenges in the General Case can be understood better.

We consider the Type I Domain here, and then handle the Type II Domain by an inversion process discussed in Section~\ref{invlace} 
Let $L \geq 1$ be an integer. We define
\begin{align*} A_0(x) & := \, \frac{x}{(x^2+1)^{3/2}}\,, && \  B_0(x) := \, \frac{1}{(x^2+1)^{3/2}}\,, \\
A_L(x) & := \, \frac{d^L}{dx^L} \frac{x}{(x^2+1)^{3/2}}\,, &&  B_L(x) := \, \frac{d^L}{dx^L} \frac{1}{(x^2+1)^{3/2}}\,. \\
\end{align*}
It is easy to see that there are polynomials $P_L$ of degree $L+1$ and $Q_L$ of degree $L$ such
that
$$A_L(x) = \frac{P_L(x)}{(1+x^2)^{(2L+3)/2}} = \frac{(1+x^2)P_{L-1}^\prime(x) - (2L+1)xP_{L-1}(x)}{(1+x^2)^{(2L+3)/2}}$$
and
$$B_L(x) = \frac{Q_L(x)}{(1+x^2)^{(2L+3)/2}} = \frac{(1+x^2)Q_{L-1}^\prime(x) - (2L+1)xQ_{L-1}(x)}{(1+x^2)^{(2L+3)/2}}\,.$$
As the degree $L+1$ of the polynomial $P_L$ and the degree $L$ of the polynomial $Q_L$ are
less than $(2L+3)/2$ we have
\begin{equation} \label{F16} \lim_{x \rightarrow \pm \infty}{\frac{P_L(x)}{(1+x^2)^{(2L+3)/2}}} = 
\lim_{x \rightarrow \pm \infty}{\frac{Q_L(x)}{(1+x^2)^{(2L+3)/2}}} = 0.
\end{equation}

\begin{prop}\label{PLPINTERLACE} Let $1 \leq L$ be an integer.

\medskip

\noindent (a) The polynomial $P_L$ of degree $L+1$ has $L+1$ distinct real zeros,
and the zeros of $P_L$ and $P_{L-1}$ strictly interlace.

\medskip

\noindent (b) The polynomial $Q_L$ of degree $L$ has $L$ distinct real zeros,
and the zeros of $Q_L$ and $Q_{L-1}$ strictly interlace.

\end{prop}

\begin{proof} We prove only (a), since the proof of (b) is identical.
The proof is a simple induction on $L$. One can check that the statements are true for $L=1$.
Assume that statement (a) is true for some integer $L\ge 1$. Note that the zeros of $P_L$ and $A_L$ are the same.
Let us denote the distinct real zeros of $P_L$ and $A_L$ by $x_1 < x_2 < \cdots < x_{L+1}$. The limit relations (\ref{F16}) 
and Rolle's theorem imply that $A_{L+1} = A^\prime_L$ has zeros $y_{j+1}\in (x_j,x_{j+1})$, $j=0,1,\ldots,L+1$, where $x_0 = -\infty$ and $x_{L+2} = \infty$.
However, the zeros of $A_{L+1}$ are the same as the zeros of $P_{L+1}$, which means that statement (a) is true for $L+1$.
\end{proof}

\begin{prop}\label{PPRIMEINTERLACE} Let $0 \leq L$ be an integer.

\medskip

\noindent (a) The zeros of $P_L$ and $P_L^\prime$ strictly interlace.

\medskip

\noindent (b) The zeros of $Q_L$ and $Q_L^\prime$  strictly interlace.

\end{prop}

\begin{proof}  We prove (b) since the proof of (a) is similiar.  The polynomial $Q_L$ is degree $L$ and has $L$ distinct zeros $x_1 < x_2 < \cdots < x_L$.
So Rolle's Theorem implies that $Q_L^\prime$ has zeros $y_j \in (x_j,x_{j+1})$, $j=1,2,\ldots,L-1$.
But the degree of $Q_L^\prime$ is $L-1$, so $Q_L^\prime$ has $L-1$ distinct zeros, and the zeros of $Q_L$ and $Q_L^\prime$ strictly interlace.
\end{proof}

\begin{prop}\label{PQINTERLACE}
Let $L \geq 1$ be an integer. The zeros of $A_L$ and $B_L$ strictly interlace.
Equivalently, the zeros of $P_L$ and $Q_L$ strictly interlace.
\end{prop}

\begin{proof}
Observe that
$$A_L(\beta) = \left. \frac{d^L A_0}{dx^L} \right|_{\beta} \qquad \text{and} \qquad A_0(\beta) = \beta B_0(\beta)\,,$$
so using the Leibniz formula we obtain
$$A_L(\beta) = LB_{L-1}(\beta) + \beta B_L(\beta)\,,$$
or equivalently
\begin{equation} \label{F17} 
P_L(\beta) = LQ_{L-1}(\beta)(1 + \beta^2) + \beta Q_L(\beta).
\end{equation}
The fact that $P_L$ and $Q_L$ have no common zeros follows from (\ref{F17}) because
any such common zero of $P_L$ and $Q_L$ is a common zero of $Q_L$ and $Q_{L-1}$, which is
impossible by Proposition~\ref{PPRIMEINTERLACE}.
It is somewhat more subtle to see that the zeros of $P_L$ and the zeros of $Q_L$ interlace.
Denoting the zeros of $Q_L$ by  $\beta_1 < \beta_2 < \cdots < \beta_L$, we can deduce from (\ref{F17}) that
$$P_L(\beta_{j-1})P_L(\beta_j) < 0\,, \qquad j=2,3,\ldots,L\,,$$
since by Proposition~\ref{PLPINTERLACE} we already know that $Q_{L-1}$ has exactly one zero in $(\beta_{j-1},\beta_j)$.
Hence $P_L$ has at least one zero in each of the intervals $(\beta_{j-1},\beta_j)\,, \enskip j=2,3,\ldots,L\,.$
By recalling the recursive formulas for $P_L$ and $Q_L$, it is easy to see by 
induction on $L$ that the sign of the leading coefficient $c_{L+1}$ of $P_L$ is $(-1)^L$, and the sign of the 
leading coefficient $d_{L-1}$ of $Q_{L-1}$ is $(-1)^{L-1}$.
Hence,
\begin{equation} \label{F18} \lim_{\beta \rightarrow \pm \infty}
\frac{P_L(\beta)}{\beta^{L+1}} = (-1)^L|c_{L+1}|\,.
\end{equation}
It follows from (\ref{F17}) and (\ref{F18}) that $P_L$ has a sign change, and hence at least one zero,
in $(\beta_L, \infty)$, and also $P_L$ has a sign change, and hence at least one zero, in $(-\infty, \beta_1)$.
Finally recall that the degree of $P_L$ is $L+1$ and the degree of $Q_L$ is $L$.
In conclusion, there is exactly one zero of $Q_L$ strictly between any two consecutive real zeros of $P_L$.
\end{proof}

\subsection{The Inversion Formula and Interlacing Properties of Polynomials}\label{invlace}

There is an interesting connection between the values of $\beta$ in the Type I Domains, and the values of $\alpha$ in the Type II Domains. 
This connection is not at first evident because of constraints that need to be made in order that we have convergence of the power series 
that we consider. But we can observe this now.

Let $L \geq 1$ be an integer. We define
\begin{align*} A_0(x) & := \, \frac{x}{(x^2+1)^{3/2}}\,, && \  B_0(x) := \, \frac{1}{(x^2+1)^{3/2}}\,, \\
A_L(x) & := \, \frac{d^L}{dx^L} \frac{x}{(x^2+1)^{3/2}}\,, &&  B_L(x) := \, \frac{d^L}{dx^L} \frac{1}{(x^2+1)^{3/2}}\,, \\ 
C_0(x) & := \, \frac{x^2}{(x^2+1)^{3/2}}\,, && \ D_0(x) := \, \frac{x}{(x^2+1)^{3/2}}\,, \\
C_L(x) & := \, \frac{d^L}{dx^L} \frac{x^{L+2}}{(x^2+1)^{3/2}}\, && D_L(x) := \, \frac{d^L}{dx^L} \frac{x^{L+1}}{(x^2+1)^{3/2}}\,. \\  
\end{align*}

\begin{prop}\label{INVERSION} (Inversion Formulas)
We have
$$\text{\rm sgn}(x)x^{L+1} B_L(x) = (-1)^L C_L(1/x)\,, \qquad x \in {\mathbb R} \setminus \{0\}\,, \qquad L=0,1,\ldots\,,$$

and

$$\text{\rm sgn}(x)x^{L+1} A_L(x) = (-1)^L D_L(1/x)\,, \qquad x \in {\mathbb R} \setminus \{0\}\,, \qquad L=0,1,\ldots\,.$$
\end{prop}

\begin{proof} We prove the first formula.  The second formula can be proved in a similar fashion.  

As both sides are even functions of $x \in {\mathbb R} \setminus \{0\}$ if $L$ is even, 
and both sides are odd functions of $x \in {\mathbb R} \setminus \{0\}$ if $L$ is odd, without loss of generality we may assume that $x > 0$.
For an integer $L \geq 0$ we define
\begin{equation} \notag 
\widetilde{B}_L(x) = \, x^{L+1} B_L(x)\,, \qquad C_L^*(x) = \, C_L(1/x)\,, \qquad \widetilde{C}_L(x)  = \, (-1)^L C_L^*(x)\,.
\end{equation}
Observe that for $L \geq 1$ we have
\begin{equation} \notag \begin{split} \widetilde{B}_L(x) & = \, x^{L+1} B_L(x) = x  x^L B_{L-1}^\prime(x) \cr 
& = \, x(x^L B_{L-1}^\prime(x)) + x(Lx^{L-1} B_{L-1}(x)) - Lx^L B_{L-1}(x) \cr 
& = \, x\widetilde{B}_{L-1}^{\prime}(x) - L\widetilde{B}_{L-1}(x)\,, \cr \end{split}
\end{equation}
that is, the functions $\widetilde{B}_L$ satisfy the recursion
\begin{equation} \label{F19} \widetilde{B}_0(x) = \frac{x}{(1+x^2)^{3/2}}\,, 
\end{equation}
\begin{equation} \label{F20}
\widetilde{B}_L(x) =  x\widetilde{B}_{L-1}^{\prime}(x) - L\widetilde{B}_{L-1}(x)\,, \qquad  L=1,2,\ldots\,.  
\end{equation}
Observe that
\begin{equation} \label{F21} \widetilde{C}_0(x) = \frac{x}{(x^2+1)^{3/2}}\,, 
\end{equation}
and the Leibniz formula yields
\begin{equation} \label{F22} C_L(x) = \frac{d^L}{dx^L} \frac{x^{L+1}x}{(x^2+1)^{3/2}} = xC_{L-1}^\prime(x) + LC_{L-1}(x)\,, \qquad L=1,2,\ldots\,. 
\end{equation}
Replacing $x$ by $1/x$ in (\ref{F22}) we get
$$C_L(1/x) = (1/x)C_{L-1}^\prime(1/x) + LC_{L-1}(1/x)\,, \qquad L=1,2,\ldots\,,$$
that is,
$$C_L^*(x) = -xC_{L-1}^{*\prime}(x) + LC_{L-1}^*(x)\,, \qquad L=1,2,\ldots\,,$$
and hence
\begin{equation} \label{F23}
\widetilde{C}_L(x) =  x\widetilde{C}_{L-1}^{\prime}(x) - L\widetilde{C}_{L-1}(x)\,, \qquad  L=1,2,\ldots\,.  
\end{equation}
Now observe that by (\ref{F19}), (\ref{F20}), (\ref{F21}), and (\ref{F23}) the functions $\widetilde{B}_L$ and $\widetilde{C}_L$ satisfy the same
recursion. In conclusion
$$\widetilde{B}_L(x) = \widetilde{C}_L(x)\,, \qquad L=0,1,\ldots\,,$$
and the proposition follows.
\end{proof}

The result of Proposition~\ref{INVERSION} is that the results of Section~\ref{lace} that applied to the Type I Domain can be transferred to the analogous properties in the Type II Domain.  For example, we have this proposition.  The other properties in the Type II Domain follow in a similar fashion.

\begin{prop}\label{CDNOTCOMMON}
Let $L \geq 1$ be an integer. The functions $C_L$ and $D_L$ do not have a common zero different from $0$.
\end{prop}

\begin{proof}  By the Inversion Formulas given by Proposition~\ref{INVERSION}, the map $1/x$ gives a one-to-one correspondence between the nonzero roots of $A_L$ and the nonzero roots of $D_L$, and 
between the nonzero roots of $B_L$ and the nonzero roots of $C_L$. Recall that by Proposition~\ref{PQINTERLACE} the zeros of  $A_L$ and the zeros of $B_L$ strictly interlace. Hence the  functions $C_L$ and $D_L$ do not have a common zero different from $0$.
\end{proof}

\subsection{The Boundary Case $y = \pm x$} \label{Boundary}
In a forthcoming paper we may be able to prove that $y = \pm x$ cannot be a common asymptotic direction to both of the zero sets $X=0$ and $Y=0$. 
Originally we needed this to argue that the zero set $X=Y=0$ is bounded. However, Proposition~\ref{FINITE} already implies that
in the Special Case the zero set $X=Y=0$ is bounded, as it is a finite set of points.

\subsection{Summary of the Special Case}\label{SpecialSum}

In this section we prove Proposition~\ref{DISTINCT}.

\begin{proof} Combining the results in Section~\ref{SecYTypeI} through Section~\ref{Boundary}, we conclude that in the Special Case the set of possible 
asymptotic directions of the zero set $X=0$ and the set of possible asymptotic directions of the zero set $Y=0$ 
are distinct, except possibly the lines $y = \pm x$. 
\end{proof}

Whether or not the lines $y = \pm x$ can be a common asymptotic direction to both of the zero sets 
$X=0$ and $Y=0$ remains open in this paper.   

The results above suggest, but do not prove, that the zero set of $X$ or the zero set of $Y$ in the Special Case are actually algebraic sets. That is, for example in the case of $X$, 
the question is whether there is a polynomial $P(x,y)$ such that the zero set $X=0$ and the zero set $P=0$.  This would mean, in particular, that where there are points in 
the zero set $X=0$ going to infinity with a particular asymptotic direction $\gamma$, then there, in fact, must be an algebraic curve in the zero set $X$ going to infinity 
(an infinite branch of the zero set $P=0$ with the same asymptotic direction $\gamma$.

\section{Charges in a General Position on the Plane}\label{GENERAL}

We now focus on a study the structure of the zero set of a finite point charge electric vector field ${\mathbf F}=(X,Y)$ in 
${\mathbb R}^2$ for which the $M$ point charges can be a general position.  We establish equations satisfied by the possible directions for the zero sets $X=0$ and $Y=0$ 
separately, and we show that there are only finitely many possible asymptotic directions for both of these 
zero sets, given again by the zero sets of sequences of certain polynomials.  
One result that these computations give, which does not seem to be obvious without the computations, is that the component zero sets have a finite number of asymptotic directions with a bound on the number of directions that is $O(M)$. 

In contrast to the Special Case, the polynomials in question are more complicated.  Indeed, the asymptotic directions, and the related of zeros of the polynomials, are not interlacing in general.  
Moreover, we do not even know if the set of asymptotic directions for the zero set $X=0$ and the set of asymptotic directions for $Y=0$ are (essentially) distinct.  These polynomials 
and their structure may have independent interest, so we include the details for one of the derivations and leave the others to be completed by analogy.  

Again, we call the domains $\{(x,y):|x| < |y|\}$ and $\{(x,y):|y| < |x|\}$ Type I Domain and Type II Domain, respectively.
In fact, for technical reasons, for a fixed $\delta \in (0,1/2)$, we will consider Type I Domains
$\{(x,y):|x| < (1-\delta)|y|\}$ and Type II Domains $\{(x,y):|y| < (1-\delta)|x|\}$.  

\subsection{Type I Domain and the Zero Set $X=0$}\label{GenXI}

The technical differences between the proof of the next proposition and the analogous results in the Special Case make it worthwhile giving this proof here. 
The other three related results in the General Case that follow below have similar results with various technical differences.  But these will be omitted to 
keep the article from getting even longer.

\medskip

\noindent {\bf Notation}:  In order to make equations more compact, we will use the notation $S(l,i) = \sum_{j=1}^M a_j(-x_j)^ly_j^i$ in this article.
The number $0 \leq L \leq 2M-1$ introduced by Definition~\defn{LGen} plays a crucial role in the formulation of our main results in this section. 

\begin{prop}\label{GenXTypeI} Let $\delta \in (0,1/2)$ be fixed. If $X(p_m,q_m) = 0$, where $(p_m,q_m)$ is in
$$ D_{\delta} = \left\{(x,y) \in \mathbb{R}^2: \left| \frac xy \right| \leq 1-\delta, \quad \frac{|x|+|x_j|}{|y|-|y_j|} < 1, \quad |y| > |y_j|, \quad j=1,2,\ldots,M \right\}\,,$$
$$|q_m| \rightarrow \infty \qquad \text{and} \qquad \lim_{m \rightarrow \infty}{\frac{p_m}{q_m}} = \beta\,,$$
then
\begin{equation} \notag
0 = \sum_{l=0}^L S(l,L-l) \frac {1}{l!(L-l)!} \frac{d^{L-l}}{dz^{L-l}} 
\left( \frac {d^l}{dz^l} \left( \frac{z}{(1 +z^2)^{3/2}} \right)z^{L+1} \right)\bigg |_\beta\,, 
\end{equation}
where $0 \leq L$ is defined by Definition~\ref{LGen}.
There are at most $2L+1 \leq 4M-1$ values of $\beta \neq 0$ satisfying the above equation.
\end{prop}

\begin{rem}
We could also interpret the zeroing identity by reversing the order that derivatives are taken. That is, we rescale the equation by $\beta^{L-l+2}$
and interpret the factors 
$$(L-l+2n+2)(L-l-1+2n+2) \cdots (1+2n+2)$$ 
as indicating $L-l$ derivatives have been taken, which would leave terms $\beta^{2n+2}$. But
then rescaling by $\beta^{-1}$ would give the zeroing identity
$$0 = \sum_{l=0}^L S(l,L-l) \frac {1}{l!(L-l)!} \frac {d^{l}}{dz^{l}} 
\left( \frac{d^{L-l}}{dz^{L-l}}\left( \frac{z^{L-l+2}}{(1 +z^2)^{3/2}} \right) \frac {1}{z} \right)\bigg |_\beta\,.$$
Note that this precludes using $\beta = 0$.
\end{rem}

\begin{rem} It is not clear which of these two forms is best to use, and it actually is not so obvious why they give the same zero set, except for $\beta = 0$.
\end{rem}

\rm Now we prove Proposition~\ref{GenXTypeI}.

\begin{proof}
We can consider a Type I Domain and consider the Taylor expansion for the zeroing equation for $X$. 
As with previous series expansions, we need to use both the binomial and geometric series. By the Binomial Theorem we have 
$$\frac 1{(1+t)^{3/2}} = (1-t)^{-3/2} = \sum_{n=0}^\infty \binom{-3/2}{n} t^n\,, \qquad t \in (-1,1)\,.$$
Hence, for $\displaystyle{\left| \frac{x-x_j}{y-y_j} \right| < 1}$ and $|y| > |y_j|$,
\begin{equation} \notag \begin{split}
\pm X(x,y) & = \, \sum_{j=1}^M{\frac{a_j(x-x_j)}{(y-y_j)^3} \left( 1 + \left( \displaystyle{ \frac{x-x_j}{y-y_j}} \right)^2 \right)^{-3/2}} \cr
& = \sum_{n=0}^\infty{\binom{-3/2}{n} \sum_{j=1}^M{\frac{a_j(x-x_j)^{2n+1}}{(y-y_j)^{2n+3}}}} \cr 
& = \, \sum_{n=0}^\infty{\binom{-3/2}{n} \frac{1}{y^{2n+3}} \sum_{j=1}^M{\frac{a_j(x-x_j)^{2n+1}}{(1-y_j/y)^{2n+3}}}}. \cr
\end{split} \end{equation}
This equation is like the one for $X$ from the Special Case analysis.  Indeed, now both of the zero sets $Y=0$ and $X=0$ have the same combined nature,
a need for a binomial expansion for the numerator terms $(x-x_j)^{2n}$ and a geometric powers series expansion for $(1-y_j/y)^{-(2n+3)}$.
So we carry out these expansions.

The Binomial Theorem gives that for $|y| > |y_n|$,
$$\frac{1}{(1-y_j/y)^{2n+3}} = \sum_{i=0}^\infty{\binom{-(2n+3)}{i}\left( \frac{-y_j}{y} \right)^i} = 
\sum_{i=0}^\infty \binom{i+2n+2}{2n+2}\left (\frac {y_j}y\right )^i\, .$$
So we have $X(x,y)$ equal to
\begin{equation} \label{GenF1} \begin{split}
& \sum_{n=0}^\infty{\binom{-3/2}{n} \frac{1}{y^{2n+3}} \sum_{j=1}^M{a_j\sum_{l=0}^{2n+1}\binom{2n+1}{l}(-x_j)^lx^{2n+1-l}\sum_{i=0}^\infty \binom{i+2n+2}{2n+2}\left (\frac {y_j}y\right )^i}} \cr
& = \, \sum_{n=0}^\infty{\binom{-3/2}{n} \sum_{l=0}^{2n+1}\sum_{i=0}^\infty \binom{2n+1}l \binom{i+2n+2}{2n+2} S(l,i) \frac {x^{2n+1-l}}{y^{2n+3 +i}}} \cr 
\end{split} \end{equation}
for all $(x,y) \in D_\delta$. Now for our asymptotic analysis, we need to identify values of $l$ and $i$ such that the moment $S(l,i) = \sum_{j=1}^M a_j(-x_j)^ly_j^i$ is zero and not zero.

To eliminate the terms that are zero, we take the moments with smallest $L = l+i$ for which at least one of these moments is not zero, 
say with $l= l_0$ and $i=i_0$.  So we have 
$$\sum_{n=0}^\infty{\binom{-3/2}{n} \sum_{l=0}^{2n+1}\sum_{i=0}^\infty \binom{2n+1}l \binom{i+2n+2}{2n+2} S(l,i)
\frac{x^{2n+1}}{y^{2n+1}}\frac {x^{-l}}{y^{-l}}\frac 1{y^{l+i+2}}} $$
equal to $0$ for all $(x,y) \in D_\delta$, where the sums over $l$ and $i$ are over all pairs $(l,i)$ with $l+i \geq L = l_0+i_0$.
The above triple sum converges absolutely in $D_\delta$ as
\begin{equation} \notag
 \sum_{n=0}^\infty{\left| \binom{-3/2}{n} \right| \sum_{l=0}^{2n+1}\sum_{i=0}^\infty \binom{2n+1}l \binom{i+2n+2}{2n+2}
|S(l,i)| 
\left| \frac{x^{2n+1}}{y^{2n+1}}\frac {x^{-l}}{y^{-l}}\frac 1{y^{l+i+2}}\right|}
\end{equation} 
is no larger than
\begin{equation}\notag \begin{split}
& \sum_{n=0}^\infty \left| \binom{-3/2}{n} \right| \sum_{l=0}^{2n+1} \sum_{i=0}^\infty \binom{2n+1}l \binom{i+2n+2}{2n+2}  \sum_{j=1}^M{|a_j||x_j|^l|y_j|^i} 
\left| \frac{x^{2n+1}}{y^{2n+1}} \frac{x^{-l}}{y^{-l}} \frac{1}{y^{l+i+2}} \right| \cr
&= \, \sum_{j=1}^M{\frac{|a_j|(|x|+|x_j|)}{\left( (|x|+|x_j)^2 + (|y| - |y_j|)^2 \right)^{3/2}}}\,, \qquad (x,y) \in D_\delta\,. \cr
\end{split} \end{equation}
Now assume that $(p_m,q_m) \in D_{\delta}$, $X(p_m,q_m) = 0$, $|q_m| \geq 1$, $m=1,2,\ldots$,  
$$\lim_{m \rightarrow \infty}{\frac{p_m}{q_m}} = \beta \in [-1+\delta,1-\delta] \qquad \text{and} \qquad \lim_{m \rightarrow \infty}{|q_m|} = \infty\,.$$ 
We plug $(p_m,q_m)$ in (\ref{GenF1}), multiply by $q_m^{L+2}$, and separate the terms in which $l+i=k=L$, in which $L+1 \leq l+i = k \leq 2n+1$, and in which $l+i = k \geq N = \max\{2n+2,L+1\}$ 
to obtain the following.  As a reminder, we are letting $S(l,i) = \sum\limits_{j=1}^M{a_j(-x_j)^ly_j^i}$.

\begin{equation}\label{GenF2} \begin{split}
0 = &  q_m^{L+2} X(p_m,q_m) \cr 
= & \sum_{n=0}^\infty \binom{-3/2}{n} \sum_{l=0}^{2n+1} \sum_{i=0}^\infty \binom{2n+1}{l} \binom{i+2n+2}{2n+2}  S(l,i) \left( \frac{p_m}{q_m} \right)^{2n+1-l} \frac{q_m^{L+2}}{q_m^{l+i+2}} \cr
= & \sum_{n=0}^\infty{\binom{-3/2}{n} \sum_{l=0}^L \binom{2n+1}{l} \binom{L-l+2n+2}{2n+2} S(l,L-l) \left( \frac{p_m}{q_m} \right)^{2n+1-l}} \cr 
+ & \sum_{n=0}^\infty{\binom{-3/2}{n} \sum_{k=L+1}^{2n+1}\sum_{l=0}^k \binom{2n+1}{l} \binom{k-l+2n+2}{2n+2} 
S(l,k-l) \left( \frac{p_m}{q_m} \right)^{2n+1-l}} \frac{1}{q_m^{k-L}} \cr
+ &  \sum_{n=0}^\infty{\binom{-3/2}{n} \sum_{k=N}^{\infty}\sum_{l=0}^{2n+1} \binom{2n+1}{l} \binom{k-l+2n+2}{2n+2} S(l,k-l) \left( \frac{p_m}{q_m} \right)^{2n+1-l}} \frac{1}{q_m^{k-L}}\,. \cr
\end{split} \end{equation}
We now show that the last two sums in (\ref{GenF2}) multiplied by $|q_m|^{1/2}$ is uniformly bounded when $(p_m,q_m) \in D_{\delta}$ and $|q_m| \geq c$ 
with a constant $c$ depending only on $\delta$, 
$$\gamma := 1 + \max_{1 \leq j \leq M}{(|x_j| + |y_j|)} \qquad \text{and} \qquad A := \sum\limits_{j=1}^M |a_j|\,.$$
We will need the estimate
\begin{equation} \label{GenF3} 
\sum_{n=0}^\infty{(n+1)^kt^n} \leq \sum_{n=0}^\infty{(n+1)(n+2) \cdots (n+k)t^n} = \frac{k!}{(1-t)^{k+1}}, \qquad t \in [0,1)\,.
\end{equation}
Observe that $\displaystyle{\frac{k-L-1/2}{2L+4} \geq k+1}$ for $k \geq L+1$, and hence with the notation $s_m = \left|q_m\right|^{1/(2L+4)}$  we have 
\begin{equation}\label{GenF4}\begin{split}
& \left| \sum_{n=0}^\infty \binom{-3/2}{n} \sum_{k=L+1}^{2n+1} \sum_{l=0}^k \binom{2n+1}{l} \binom{k-l+2n+2}{2n+2} 
S(l,k-l) ( \frac{p_m}{q_m} )^{2n+1-l} \frac{|q_m|^{1/2}}{q_m^{k-L}} \right| \cr
\leq  & \sum_{n=0}^\infty (n+1) \sum_{k=L+1}^{2n+1} \sum_{l=0}^k{\frac{(2n+1)^{l}}{l!} \frac{(k-l+2n+2)^{k-l}}{(k-l)!}} \frac{A\gamma(1-\delta)^{2n-1-l}}{|q_m|^{k-L-1/2}} \cr
\leq  & \sum_{n=0}^\infty (n+1) \sum_{k=L+1}^{2n+1} \sum_{l=0}^k{\frac{(4n+3)^k}{l!(k-l)!} A\gamma^k (1-\delta)^{2n+1-l}} \frac{1}{s_m^{k+1}} \cr
\leq  & \sum_{n=0}^\infty (n+1) \sum_{k=L+1}^{2n+1} \sum_{l=0}^k{\frac{1}{k!} \binom{k}{l} (4n+3)^k A\gamma^k (1-\delta)^{2n+1-l}} \frac{1}{s_m^{k+1}} \cr
\leq  & \sum_{n=0}^\infty (n+1) \sum_{k=L+1}^{2n+1} \frac{2^k}{k!} (4n+3)^k A\gamma^k (1-\delta)^{2n+1-l} \frac{1}{s_m^{k+1}}\,. \cr
\end{split} \end{equation}
Now we break the sum in the last line of (\ref{GenF3}) for $L+1 \leq k \leq n$ and $n+1 \leq k \leq 2n+1$. Recalling (\ref{GenF3}),
we have  
\begin{equation} \label{GenF5} \begin{split}
& \, \sum_{n=0}^\infty (n+1) \sum_{k=L+1}^{n} \frac{2^k}{k!} (4n+3)^k A\gamma^k (1-\delta)^{2n+1-l} \frac{1}{s_m^{k+1}} \cr 
\leq & \sum_{n=0}^\infty \sum_{k=L+1}^{n} \frac{2^k}{k!} (4(n+1))^{k+2} A\gamma^k (1-\delta)^{n+1} \frac{1}{s_m^{k+1}} \cr
\leq & \sum_{n=0}^\infty \sum_{k=L+1}^{\infty} A \frac{(8\gamma(n+1))^{k+2}}{k!} (1-\delta)^{n+1} \frac{1}{s_m^{k+1}} \cr  
= & \sum_{k=L+1}^\infty \sum_{n=0}^{\infty} A \frac{(8\gamma(n+1))^{k+1}}{k!} (1-\delta)^{n+1} \frac{1}{s_m^{k+1}} \cr
= & \sum_{k=L+1}^\infty \sum_{n=0}^{\infty} \frac{A}{k!} \left( \frac{8\gamma(n+1)}{s_m} \right)^{k+2} (1-\delta)^{n+1} \cr
\leq & \sum_{k=L+1}^\infty \sum_{n=0}^{\infty} \frac{A}{(k+1)!} (n+1)^{k+1} (1-\delta)^{n+1} \left( \frac{16\gamma}{s_m} \right)^{k+1} \cr
\leq & \sum_{k=L+1}^\infty A \left( \frac{1}{\delta} \right)^{k+2} \left( \frac{16\gamma}{s_m} \right)^{k+2} \cr
\leq & \sum_{k=L+1}^\infty A \left( \frac{16\gamma}{s_m\delta} \right)^{k+2} \leq B_1\,,\cr
\end{split} \end{equation}
with a constant $B_1$ depending only on $\delta$, $A$, and $\gamma$ if $s_m = \left|q_m\right|^{1/(2L+4)} \geq 32\gamma/\delta$. 
 
Similarly
\begin{equation}\label{GenF6} \begin{split}
& \sum_{n=0}^\infty (n+1) \sum_{k=n+1}^{2n+1} \frac{2^k}{k!} (4n+3)^k A\gamma^k (1-\delta)^{2n+1-l} \frac{1}{s_m^{k+1}} \cr 
\leq & \sum_{n=0}^\infty \sum_{k=n+1}^{2n+1} \frac{2^k}{k!} (4(n+1))^{k+1} A\gamma^k \frac{1}{s_m^{k+1}} 
\leq \sum_{n=0}^\infty \sum_{k=n+1}^{2n+1} A \frac{(8\gamma(n+1))^{k+2}}{k!} \frac{1}{s_m^{k+2}} \cr  
\leq & \sum_{n=0}^\infty \sum_{k=n+1}^{2n+1} A \left( \frac{16e\gamma(n+1)}{ks_m} \right)^{k+1} 
\leq \sum_{n=0}^\infty \sum_{k=n+1}^{2n+1} A \left( \frac{16e\gamma(n+1)}{{(n+1)s_m}} \right)^{k+1} \cr 
\leq & \sum_{n=0}^\infty \sum_{k=n+1}^{2n+1} A \left( \frac{32e\gamma}{s_m} \right)^{k+1}2^{-n} \leq B_2 \cr
\end{split} \end{equation}
with a constant $B_2$ depending only on $\delta$, $A$, and $\gamma$ if $s_m = \left|q_m\right|^{1/(2L+4)} \geq 64e\gamma$.
Combining (\ref{GenF4}), (\ref{GenF5}), and (\ref{GenF6}) completes the proof that the last but one sum  
\begin{equation} \notag 
\sum_{n=0}^\infty{\binom{-3/2}{n} \sum_{k=L+1}^{2n+1}\sum_{l=0}^k \binom{2n+1}{l} \binom{k-l+2n+2}{2n+2} 
S(l,k-l) \left( \frac{p_m}{q_m} \right)^{2n+1-l}} \frac{1}{q_m^{k-L}} 
\end{equation}
in (\ref{GenF2}) converges to $0$ when $\displaystyle{\lim_{m \rightarrow \infty}{|q_m|}} = \infty$.

With $N = \max\{2n+2,L+1\}$ and $s_m = \left|q_m\right|^{1/(2L+4)}$ we also have
\begin{equation} \label{GenF7} \begin{split}
& \, \left| \sum_{n=0}^\infty{\binom{-3/2}{n} \sum_{k=N}^\infty \sum_{l=0}^{2n+1} \binom{2n+1}{l} \binom{k-l+2n+2}{2n+2} 
S(l,k-l) \left( \frac{p_m}{q_m} \right)^{2n+1-l}} \frac{|q_m|^{1/2}}{q_m^{k-L}} \right| \cr
\leq & \,\sum_{n=0}^\infty (n+1) \sum_{k=N}^\infty \sum_{l=0}^{2n+1} \binom{2n+1}{l} \binom{k-l+2n+2}{2n+2}  A\gamma^k \frac{1}{|q_m|^{k-L-1/2}} \cr
\leq & \, \sum_{n=0}^\infty (n+1) \sum_{k=N}^{\infty} \sum_{l=0}^{2n+1} \binom{2n+1}{l} \frac{(k+2n+2)^{2n+2}}{(2n+2)!} A\gamma^k \frac{1}{s_m^{k+1}} \cr
\leq & \, \sum_{n=0}^\infty (n+1) \sum_{k=N}^{\infty} 2^{2n+1} \frac{(k+2n+2)^{2n+2}}{(2n+2)!} A\gamma^k  \frac{1}{s_m^{k+1}} \cr
\leq & \, \sum_{n=0}^\infty \sum_{k=N}^{\infty} (n+1) \left( \frac{2e(2n+2+k)}{2n+2} \right)^{2n+2} A\gamma^k  \frac{1}{s_m^{k+1}} \cr
\leq & \, \sum_{n=0}^\infty \sum_{k=N}^{\infty} (n+1) (2e)^{2n+2}\left( 1 + \frac{k}{2n+2} \right)^{2n+2} A\gamma^k  \frac{1}{s_m^{k+1}} \cr
\leq & \, \sum_{n=0}^\infty \sum_{k=N}^{\infty} (n+1) (2e)^{2n+2} e^k A\gamma^k  \frac{1}{s_m^{k+1}} 
\leq \sum_{n=0}^\infty \sum_{k=N}^{\infty} (n+1) (2e)^{2n+2} A \left( \frac{e\gamma}{s_m} \right)^k \frac{1}{s_m} \cr  
\leq & \, \sum_{n=0}^\infty \sum_{k=N}^{\infty} (n+1) A \left( \frac{4e^2\gamma}{s_m} \right)^k (2e)^{-(2n+2)}\frac{1}{s_m} \leq B_3 \cr
\end{split} \end{equation}
with a constant $B_3$ depending only on $\delta$, $A$, and $\gamma$ if $s_m = \left|q_m\right|^{1/(2L+4)} \geq 8e^2\gamma$. It follows from (\ref{GenF7}) that the last sum 
\begin{equation} \notag
\sum_{n=0}^\infty{\binom{-3/2}{n} \sum_{k=2n+2}^{\infty} \sum_{l=0}^{2n+1} \binom{2n+1}{l} \binom{k-l+2n+2}{2n+2} 
S(l,k-l) ( \frac{p_m}{q_m} )^{2n+1-l}} \frac{1}{q_m^{k-L}} 
\end{equation}
in (\ref{GenF2}) converges to $0$ when $q_m$ tends to $\infty$.

Now we are ready to take the limit in (\ref{GenF2}) when $\displaystyle{\lim_{m \rightarrow \infty}{|q_m|} = \infty}$ and $\displaystyle{\lim_{m \rightarrow \infty}{p_m/q_m} = \beta\,.}$ 
We have seen that the last two sums in (\ref{GenF2}) tend to $0$ as $\displaystyle{\lim_{m \rightarrow \infty}{|q_m|} = \infty}$. So we have 
\begin{equation} \notag \begin{split}
0 & = \, \lim_{|q_m| \rightarrow \infty} \sum_{n=0}^\infty{\binom{-3/2}{n} \sum_{l=0}^L \binom{2n+1}{l} \binom{L-l+2n+2}{2n+2} S(l,L-l)
\left( \frac{p_m}{q_m} \right)^{2n+1-l}} \cr  
& = \, \sum_{l=0}^L S(l,L-l) \beta^{-l} \sum_{n=0}^\infty \binom{-3/2}{n} \binom{2n+1}{l} \binom{L-l+2n+2}{2n+2} \beta^{2n+1}\,. \cr
\end{split} \end{equation}
Observe that the right-hand side is a not identically zero Laurent series in $0 \neq \beta \in (-1,1)$.  
We would like to show this zeroing formula is equivalent to a formula involving derivatives in such a manner that it can be seen that there are only a finite number of solutions.
We can do this by observing that
\begin{equation} \notag  
\sum_{n=0}^\infty{\binom{-3/2}{n} \binom{2n+1}{l} \binom{L-l+2n+2}{2n+2} \beta^{2n+1}}
\end{equation} 
is the same as
\begin{equation} \notag 
\frac {1}{l!(L-l)!} \sum_{n=0}^\infty \binom{-3/2}{n} (2n+1) \cdots (2n +1 -l+1) \cdot (L-l+2n+2)\cdots (1+2n+2) \beta^{2n+1}\,.
\end{equation} 
Interpret the factors 
$$(2n+1)(2n+1-1) \cdots (2n+1-l+1)$$ 
as indicating that the $l$th derivatives of $\displaystyle{\frac z{(1+z^2)^{3/2}}}$ have been taken, and evaluated at $\beta$, which would leave terms with a factor $\beta^{2n - l +1}$.
But then rescaling by $\beta^{L+1}$, would give this term to be
$$\frac 1{z} \frac {d^{L-l}}{dz^{L-l}} \left (\frac {d^l}{dz^l}\left (\frac {z}{(1 +z^2)^{3/2}}\right )z^{L+1}\right )\bigg |_\beta\,.$$
So our zeroing identity is
\begin{equation} \notag
\sum_{l=0}^L S(l,L-l) \frac {1}{l!(L-l)!} \frac{d^{L-l}}{dz^{L-l}} 
\left( \frac {d^l}{dz^l} \left( \frac{z}{(1+z^2)^{3/2}} \right)z^{L+1} \right)\bigg |_\beta = 0\,.
\end{equation}
This shows that the zeroing equation is of the form $P(\beta)/(1 +\beta^2)^{(2L+3)/2}$ for a suitable not identically zero polynomial $P$ of degree at most $2L+1 \leq 4M-1$.
So the zeroing identity only gives zeros in the roots of $P$. As $P$ is a polynomial of degree at most $2L+1 \leq 4M-1$, there are at most $2L+1 \leq 4M-1$ 
possible asymptotic directions for the zero set $X=0$ in the Type I Domain.
\end{proof}

\subsection{Type I Domain and the Zero Set $Y=0$}\label{GenYI}

There is more symmetry in our arguments with charges for the General Case. When we switch from the zero set $X= 0$ to the zero set $Y=0$ in the Type I Domain, 
it just decreases the power in the binomial terms by $1$ and decreases the powers in the geometric terms by $1$. With the necessary adjustment of the powers,
in the Type I Domain for the zero set $Y=0$ we have the zeroing identity 
$$0 = \sum_{l=0}^L S(l,L-l) \beta^{-l} \sum_{n=0}^\infty \binom{-3/2}{n} \binom{2n}l \binom{L-l+2n+1}{2n+1} \beta^{2n}\,,$$
and so the following result holds. 

\begin{prop}\label{GenYTypeI} Let $\delta \in (0,1/2)$ be fixed. If $Y(p_m,q_m) = 0$, where $(p_m,q_m)$ is in 
$$ D_{\delta} = \left\{(x,y) \in {\mathbb R}^2: \left| \frac xy \right| \leq 1-\delta, \quad \frac{|x|+|x_j|}{|y|-|y_j|} < 1, \quad |y| > |y_j|, \quad j=1,2,\ldots M,\right\}$$
$$|q_m| \rightarrow \infty \qquad \text{and} \qquad \lim_{m \rightarrow \infty}{\frac{p_m}{q_m}} = \beta\,,$$
then
\begin{equation} \notag
0 = \sum_{l=0}^L S(l,L-l) \frac{1}{l!(L-l)!} 
\frac{d^{L-l}}{dz^{L-l}} \left( \frac{d^l}{dz^l} \left( \frac{1}{(1 +z^2)^{3/2}} \right)z^{L+1} \right)\bigg |_\beta\,,
\end{equation}
where $0 \leq L$ is defined by Definition~\ref{LGen}.
There are at most $2L+1 \leq 4M-1$ values of $\beta \neq 0$ satisfying the above equation.
\end{prop}

\subsection{Common Asymptotic Directions of the Zero Sets $X$ and $Y$ in Type I Domain}\label{COMMON}

We can ask: is there a $\beta$ which gives common asymptotic directions for both of the zero sets $X=0$ and $Y=0$ in the Type I Domain?    
Such a $\beta \in [-1+\delta,1-\delta]$ must satisfy both of the equations 
$$0 = \sum_{l=0}^L \Gamma(l) 
\frac{d^{L-l}}{dz^{L-l}} \left( \frac{d^l}{dz^l} \left( \frac {z}{(1 +z^2)^{3/2}} \right)z^{L+1} \right)\bigg |_\beta$$ 
and
$$0 = \sum_{l=0}^L \Gamma(l)
\frac{d^{L-l}}{dz^{L-l}} \left( \frac{d^l}{dz^l} \left( \frac{1}{(1 +z^2)^{3/2}} \right)z^{L+1} \right)\bigg |_\beta\,.$$ 
with some real numbers $\Gamma(l)$ derived above.

We have come close to showing that this cannot hold in the Special Case where the point charges are on the positive $x$-axis. For example, 
we showed that in the Special Case the possible asymptotic directions for the zero set $X=0$ and the possible asymptotic directions for the
zero set $Y=$ are strictly interlacing in the Type I and Type II Domain.     

Examples in the General Case show that sometimes the zero sets $X=0$ and $Y=0$ have asymptotic directions that are not interlacing.  
But we can ask: do the component zero sets have distinct asymptotic directions in the Type I Domain?  If so, outside some large disk $X$ and $Y$ 
have no common zeros in the Type I Domain.  That is, we can try to use specific formulas for the asymptotic directions to show that   
${\mathbf F}=(X,Y)$ has no zeros far from the origin in the Type I Domain.

\medskip

\subsection{Type II Domain and $X=0$}\label{GenXII}

The result below can be obtained from Proposition~\ref{GenXTypeI} by reversing the roles of $x$ and $y$ and $x_j$ and $y_j$.

\begin{prop}\label{GenXTypeII} Let $\delta \in (0,1/2)$ be fixed. If $X(p_m,q_m) = 0$, where $(p_m,q_m)$ is in
$$H_{\delta} = \left\{(x,y) \in {\mathbb R}^2: \left| \frac yx \right| \leq 1-\delta, \quad \frac{|y|+|y_j|}{|x|-|x_j|} < 1, \quad |x| > |x_j|, \quad j=1,2,\ldots,M \right\}$$
$$|p_m| \rightarrow \infty \qquad \text{and} \qquad \lim_{m \rightarrow \infty}{\frac{q_m}{p_m}} = \alpha\,,$$
then
\begin{equation} \notag
0 = \sum_{l=0}^L S(L-l,l) \frac {1}{l!(L-l)!} 
\frac{d^{L-l}}{dz^{L-l}} \left( \frac{d^l}{dz^l} \left( \frac{1}{(1 +z^2)^{3/2}} \right)z^{L+1} \right)\bigg |_\alpha\,.
\end{equation}
There are at most $2L+1 \leq 4M-1$ values of $\alpha \neq 0$ satisfying the above equation.
\end{prop}

\subsection{Type II Domain and the Zero Set $Y=0$}

The result in this section can be obtained from Proposition~\ref{GenYTypeI} by reversing the roles of $x$ and $y$ and $x_j$ and $y_j$.

\begin{prop}\label{GenYTypeII} Let $\delta \in (0,1/2)$ be fixed. If $Y(p_m,q_m) = 0$, where $(p_m,q_m)$ is in
$$H_{\delta} = \left\{(x,y) \in {\mathbb R}^2: \left| \frac yx \right| \leq 1-\delta, \quad \frac{|y|+|y_j|}{|x|-|x_j|} < 1, \quad |x| > |x_j|, \quad j=1,2,\ldots,M \right\}$$
and
$$|p_m| \rightarrow \infty \qquad \text{and} \qquad \lim_{m \rightarrow \infty}{\frac{q_m}{p_m}} = \alpha\,,$$
then
\begin{equation} \notag
0 = \sum_{l=0}^L S(L-l,l) \frac {1}{l!(L-l)!} 
\frac{d^{L-l}}{dz^{L-l}} \left( \frac{d^l}{dz^l} \left( \frac{z}{(1 +z^2)^{3/2}} \right)z^{L+1} \right)\bigg |_\alpha\,.
\end{equation}
There are at most $2L+1 \leq 4M-1$ values of $\alpha \neq 0$ satisfying the above equation.
\end{prop}

\subsection{Common Asymptotic Directions of the Components $X$ and $Y$ in Type II Domain}\label{SAME}

We can ask, as we did above in Section~\ref{COMMON}, if there are no $\alpha$ which is a common asymptotic direction for both of the zero sets $X=0$ and $Y=0$ 
in the Type II Domain. Such an $\alpha \in [-1+\delta,1-\delta]$ must satisfy both of the equations
$$0 = \sum_{l=0}^L \Delta(l) \frac{d^{L-l}}{dz^{L-l}} \left( \frac{d^l}{dz^l} \left( \frac{1}{(1 +z^2)^{3/2}} \right)z^{L+1} \right)\bigg |_\alpha$$
and
$$0 = \sum_{l=0}^L \Delta(l) \frac{d^{L-l}}{dz^{L-l}} \left( \frac{d^l}{dz^l} \left( \frac{z}{(1 +z^2)^{3/2}} \right)z^{L+1} \right)\bigg |_\alpha$$
with particular real numbers $\Delta(l)$ derived above.

See Section~\ref{COMMON} for additional remarks on why this is an important issue to resolve.

\subsection{Proof of Proposition~\ref{squaring}}  

\begin{proof} Take $R$ and put the terms over the common denominator $D^{1/2}$ where $D =B_1B_2 \cdots B_M$.  This expresses
$R = D^{-1/2} \sum_{j=1}^M{A_j D_j^{1/2}}$ where each $D_j$ is $D$ with $B_j$ divided out.  We see then that the zero set $R = 0$ is a subset of where
$S = \sum_{j=1}^M {A_jD_j^{1/2}} = 0$. Let $\Sigma$ be the collection of the $2^M$ functions $\sigma: \{1,2,\ldots,M\} \rightarrow \{-1,1\}$.
Now we take the product
$$P = \prod_{\sigma \in \Sigma}{\Bigg( \sum_{j=1}^M{\sigma(j)A_j D_j^{1/2}} \Bigg)}\,.$$
Here the product has $2^M$ factors. We show that $P$ is a nonidentically $0$ polynomial (of degree at most $3M \cdot 2^M$). Let
$${\mathcal P}(\chi_1,\chi_2,\ldots,\chi_M) = \prod_{\sigma \in \Sigma}{\bigg( \sum_{j=1}^M{\sigma(j) \chi_j} \Bigg)}\,.$$
Observe that for every fixed $k \in \{1,2,\ldots,M\}$ the product ${\mathcal P}$ is a polynomial in $\chi_k$ and its value remains the same if we replace
$\chi_k$ by $-\chi_k$, as the factors of the product remain the same, they are only permuted.
Hence ${\mathcal P}$ is a polynomial in $\chi_k^2$ for every fixed $k \in \{1,2,\ldots,M\}$. Now apply this with $\chi_j = A_jD_j^{1/2}$ to conclude that $P$
is a polynomial, indeed.
To finish the proof observe that $\sum_{j=1}^M{A_jD_j^{1/2}}$ is a factor of this product defining $P$, namely the sum when $\sigma(j) =1$
for each $j$. So the zero set $R = 0$ in ${\mathbb R}^n$ is a subset of the real algebraic variety given by $P = 0$ in ${\mathbb R}^n$.
\end{proof}

\subsection{Summary of the General Case}\label{GenSum}

In the General Case, how do we show that there can only be a finite number of zeros in a given bounded disk?   This would mean showing there cannot be a curve of zeros anywhere.   Indeed, if we consider $S = X^2 + Y^2$, then the zero set $S = 0$ is an analytic variety in the plane and hence, in this case, is a locally finite collection of points 
and curves. Suppose that we can use some of the above analysis, and extensions of it, to show that there are no points in this zero set outside a large disk.  If so, then all we have to do is show that there are no curves where $X$ and $Y$ are zero.

So it is important for us to observe that the article by Abanov, Hayford, Khavinson, and Teodorescu~\cite{AHKT} gives a short proposed proof in Proposition 4.1 that there 
are no curves in the common zero set.  

\begin{conj} The zero set of a two-dimensional finite electric vector field is bounded. 
\end{conj}

\noindent The result then would be a solution of the following

\begin{conj} The zero set of a two-dimensional finite electric vector field consists of only finitely many points.
\end{conj}
\smallskip
\noindent {\bf Acknowledgments}: 

We thank Bruce Reznick for showing us the product method for eliminating square roots in our vector field component equations in Propositionn~\ref{squaring} 
and Section~\ref{Curves}, and we thank Vilmos Totik for helpful remarks about the interlacing of zeros of in our formulas for asymptotic directions in Section~\ref{Asymptotes}. 

Finally, Dr.~Rebecca Rosenblatt thanks the National Science Foundation for supporting some of this work via Independent Research/Development release time. 
Please note that the opinions, findings, conclusions, and recommendations expressed here are the authors' and do not necessarily reflect the views of the National Science Foundation.




\bibliographystyle{amsplain}

\begin{thebibliography}{99}

\bibitem{AHKT} A. Abanov, N. Hayford, D. Khavinson, and R. Teodorescu, 
{\em Around a theorem of F. Dyson and A. Lenard: energy equilibria for point charge distributions in classical electrostatics}, preprint, 13 pages.  

\bibitem{GNS} A. Gabrielov, D. Novikov, and B. Shapiro, 
{\em Mystery of point charges}, Proc. London Math. Soc. (3) 95 (2007) 443-472.

\bibitem{Gibson} C. G. Gibson, 
{\em Elementary geometry of algebraic curves: an undergraduate introduction}. Cambridge Univ. Press, Cambridge, 1998.

\bibitem{K} K. Killian, {\em A remark on Maxwell's conjecture for planar charges}, Complex Var. Elliptic Equ. 54 (2009), no. 12, 1073-1078.

\bibitem{WI} H. Whitney, {\em Tangents to an analytic variety}, Ann. Math. 81 (1965) 496-549.

\bibitem{WII} H. Whitney, 
{\em Local properties of analytic varieties}, Differential and Combinatorial Topology, Ed. S. Cairns, Princeton Univ. Press, Princeton, New Jersey, 1965.

\bibitem{WB}  H. Whitney and F. Bruhat, 
{\em Quelques propri\'et\'es fondamentales des ensembles analytiques r\'eels}, Comment. Math. Helv. 33 (1959) 132-160.

\bibitem{WIII} H. Whitney, 
{\em Elementary structure of real algebraic varieties}, Ann.  Math. (3) 66 (1957), 545-556.

\bibitem{Jan} A. I. Yanushauskas, 
{\em On the zeros of the gradient of a harmonic function}, Dokl. Akad. Nauk SSSR 158 (1964) 547-549.  
\end{thebibliography}

\bigskip

\noindent {\bf Authors}:
\medskip

\noindent Tam\'as Erd\'elyi, Department of Mathematics, Texas A\&M University, College Station, Texas,
terdelyi@math.tamu.edu
\medskip

\noindent Joseph Rosenblatt, Department of Mathematics, University of Illinois at Urbana-Champaign, Urbana, Illinois,
rosnbltt@illinois.edu
\medskip

\noindent Rebecca Rosenblatt, Division of Equity for Excellence in STEM, National Science Foundation, Alexandria,  Virginia, 
rrosenbl@nsf.gov
\bigskip

\noindent {\bf AMS Subject Classifications}: $41A58, 14H50, 26B12$
\bigskip

\noindent {\bf Alt Text}:  This article contains text in English with numerous mathematical formulas that accompanying the text.  

\noindent {\bf Received}: June 15, 2024

\end{document}